\documentclass[12pt,twoside,a4paper]{article}
\usepackage[english]{babel}
\usepackage[utf8]{inputenc}
\usepackage{amsmath}
\usepackage{amsfonts}
\usepackage[pagebackref]{hyperref}
\usepackage{amsthm}
\usepackage{amssymb}
\usepackage{bm}
\usepackage{tikz-cd}
\usepackage{xr}
\usepackage[nameinlink]{cleveref}

\usepackage{enumitem}
\usepackage{authblk}
\usepackage[title]{appendix}
\usepackage{chngcntr}
\usepackage{apptools}
\usepackage{mathtools}
\usepackage{appendix}

\usepackage{bookmark}
\usepackage{appendix}
\definecolor{crimsonglory}{rgb}{0.75, 0.0, 0.2}
\definecolor{darkpowderblue}{rgb}{0.0, 0.2, 0.6}
\hypersetup{
	colorlinks   = true, 
	urlcolor     = darkpowderblue, 
	linkcolor    = darkpowderblue, 
	citecolor   = crimsonglory 
}

\newtheorem*{rmk}{Remark}

\newtheorem*{rmks}{Remarks}

\newtheorem{claim}{Claim}

\newtheorem{theorem}{Theorem}
\newtheorem{lemma}{Lemma}
\newtheorem{defn}{Definition}
\newtheorem{cor}{Corollary}
\newtheorem*{thm*}{Theorem}
\newtheorem{conj}{Conjecture}

\newtheorem{conv}{Convention}
\newtheorem*{notation}{Notation}
\usepackage{theoremref}

\DeclareMathOperator{\End}{End}
\DeclareMathOperator{\GL}{GL}

\DeclareMathOperator{\spec}{Spec}

\DeclareMathOperator{\aut}{Aut}

\DeclareMathOperator{\lie}{Lie}
\DeclareMathOperator{\disc}{disc}
\DeclareMathOperator{\atyp}{atyp}
\DeclareMathOperator{\codim}{codim}

\newcommand{\Sum}[2]{\displaystyle\sum_{#1}^{#2}}

\newcommand{\Bigsum}[2]{\displaystyle\bigoplus_{#1}^{#2}}

\newcommand{\Z}{\mathbb{Z}}
\newcommand{\N}{\mathbb{N}}
\newcommand{\C}{\mathbb{C}}
\newcommand{\Q}{\mathbb{Q}}
\newcommand{\R}{\mathbb{R}}
\newcommand{\F}{\mathbb{F}}

\usepackage[OT2,T1]{fontenc}
\DeclareSymbolFont{cyrletters}{OT2}{wncyr}{m}{n}
\DeclareMathSymbol{\Sha}{\mathalpha}{cyrletters}{"58}

\usepackage{tocloft}

\numberwithin{conv}{section}
\numberwithin{claim}{section}
\numberwithin{theorem}{section}
\numberwithin{defn}{section}
\numberwithin{prop}{section}
\numberwithin{cor}{section}
\numberwithin{lemma}{section}

\usepackage{enumitem}

\newcommand\coolrightbrace[2]{%
	\left.\vphantom{\begin{matrix} #1 \end{matrix}}\right\}#2}

	\title{Some cases of the Zilber-Pink conjecture for curves in $\mathcal{A}_g$}
\author{Georgios Papas}
\date{}

\begin{document}

	\maketitle
	
	\begin{abstract} Following our work in \cite{papas2022height}, we extend the height bounds established by Y. Andr\'e in his seminal research monograph \cite{andre1989g} for $1$-parameter families of abelian varieties defined over number fields. In our exposition we no longer assume that the family acquires completely multiplicative reduction at some point, as in Andr\'e's original result.
		
		As a corollary of these height bounds, we obtain unconditional results of Zilber-Pink-type for curves in $\mathcal{A}_g$, building upon recent results of C. Daw and M. Orr.
	\end{abstract}

\section{Introduction}

The main inspiration of our work in \cite{papas2022height} and here is the following height bound of Andr\'e:
\begin{theorem}[\cite{andre1989g}, $X.1$]\label{andresbound} Let $S'$ be a smooth connected curve over a number field $K$, let $s_0\in S'(K)$, and let $S=S'\backslash \{s_0\}$. 
	
	Consider $f:A\rightarrow S$ an abelian scheme also defined over the number field $K$. Assume that the geometric generic fiber $A_{\eta}$ is a simple abelian variety of odd dimension $g>1$ and that the connected component of the fiber at $s_0$ of the N\'eron model of $A$ over $S'$ is a $g$-dimensional torus.
	
	Let $h$ denote a Weil height on $S'$ and consider the set \begin{center}
		$\Sha(S):= \{s\in S(\bar{\Q}):\End^{0}(A_s)\not\hookrightarrow M_g(\Q) \}$.
	\end{center}Then there exist effectively computable constants $c_1,$ $c_2>0$ such that for all $s\in \Sha(S)$ \begin{center}
		$h(s)\leq c_1 [K(s):K]^{c_2}$.
	\end{center}
\end{theorem}

This theorem of Andr\'e has an immediate Hodge theoretic interpretation. Namely, it confirms that the height of the points in $S(\bar{\Q})$ whose corresponding fibers carry more Hodge endomorphisms than the generic fiber of the abelian family considered in this case, and furthermore satisfy the condition that $\End^{0}(A_s)\not\hookrightarrow M_g(\Q)$, is controlled by the degree of their definition. This is in agreement with the general expectation raised by the Zilber-Pink conjecture that such points are in fact finite!

We note that C. Daw and M. Orr establish analogous height bounds in the case when $g$ is even, see Theorem $8.1$ of \cite{daworr}. Namely, they consider $S'$, $s_0$, $S$, and $f:A\rightarrow S$ to be as in \Cref{andresbound} but let $g$ be $>1$, i.e. there is no condition on the parity of $g$. Their theorem though still assumes that the reduction at $s_0$ of the connected N\'eron model $A'$ of $A$ over $S'$ is ``completely multiplicative'', i.e. that $A'_{s_0}$ is a torus.

In \cite{papas2022height} we establish analogous height bounds for exceptional points in certain geometric variations of Hodge structures of arbitrary odd weight. The main goal of this article is to establish height bounds analogous to the aforementioned ones for certain exceptional points in one-parameter algebraic families of abelian varieties replacing the assumption of ``completely multiplicative reduction'' by allowing semiabelian reduction instead.

All of the aforementioned results stem from Andr\'e's G-functions method. For an introduction to G-functions see \cite{dwork}.

\subsubsection{Main result}\label{section:notation}

Let $S'$ be a smooth geometrically irreducible curve defined over a number field $K$, let $s_0\in S(K)$ be a fixed closed point and let $S:=S'\backslash \{s_0\}$. Consider an abelian scheme $f:X\rightarrow S$, also defined over $K$, with $g$-dimensional fibers, where $g\geq 2$.

Let $f':X'\rightarrow S'$ be the connected N\'eron model of $X_{K(S)}$ over $S'$, i.e. the group subscheme of the N\'eron model of $X$ over $S'$ whose fiber over any point in $S'$ is the connected component of the identity of the fiber of the N\'eron model at that point. We assume that the fiber of $X'$ over $s_0$ is a semiabelian variety over $K$ whose toric part is $h$-dimensional. In other words we assume that the Chevalley decomposition of $X'_{s_0}$ is given by \begin{equation}\label{eq:semiabelian}	
	0\rightarrow T \rightarrow X'_{s_0}\rightarrow B\rightarrow 0,	
\end{equation}where $T$ is an $h$-dimensional $K$-torus and $B$ a $(g-h)$-dimensional abelian variety defined over $K$. 

We assume that the generic special Mumford-Tate group of the variation of weight $1$ Hodge structures associated to $f:X\rightarrow S$ is $Sp_{2g}(\Q)$. This is equivalent to assuming that the map $S\rightarrow \mathcal{A}_g$ induced by $f$ is Hodge generic.

\begin{defn} We call any abelian scheme $f:X\rightarrow S$ that satisfies the above with $h\geq 2$ a \textbf{G-admissible abelian scheme over $S'$}.
	
\end{defn}

Our main result is the following height bound in the spirit of Andr\'e's original result.

\begin{theorem}\label{application2}Let $S'$ be a smooth geometrically irreducible curve over the number field $K$ and let $S$ and $s_0$ be as above. Consider $f:X\rightarrow S$ a $G$-admissible abelian scheme over $S'$.
	
	We consider the set \begin{center}	$\Sha(S):= \{ P \in S(\bar{\Q}) : P \text{ is strongly exceptional } \}$.
	\end{center}Then, there exist positive constants $C_1$, $C_2$ such that $h(P)\leq C_1 [K(P):K]^{C_2}$ for all $P\in \Sha(S)$.

\end{theorem}

\begin{rmks}1. As noted earlier, our main point of divergence from Andr\'e's result and the result of Daw and Orr, is that we do not require the degeneration at $s_0$ to be completely multiplicative and, in comparison to Andr\'e's result, that we do not require any conditions about the parity of the dimension $g$.\\
	
	2. For the definition of ``strongly exceptional points''  see \Cref{section:strongexception}. The restrictions on the points we are considering, reflected in the definition of ''strong exceptionality'', are more stringent than the very elegant condition ``$\End^{0}(X_s) \not\hookrightarrow M_g(\Q)$'' of Andr\'e. 
	
	This necessity arises from two issues. On the one hand is the fact that we have fewer G-functions appearing as relative periods, at least in full generality when $h<g$. On the other hand is the increased possible complexity of the degeneration at $s_0$, which here, unlike Andr\'e's original result, does not have to be completely multiplicative.\\ 
	
	3. We note that the constants $C_i$ in the previous theorem are in fact effective. For more around the effectivity of our results see the remarks in $1.E$ of \cite{daworr3}.
\end{rmks}

\subsection{Towards the Zilber-Pink Conjecture}

The main motivation behind our pursuit of such height bounds is a recent series of papers, see \cite{daworr, daworr2, daworr3, daworr4}, by C. Daw and M. Orr, where cases of the Zilber-Pink Conjecture for curves in moduli spaces of abelian varieties are established. 

The Zilber-Pink conjecture is a far-reaching conjecture in the general context of problems of unlikely intersections. In various contexts it generalizes various previous conjectures in the field such as the Andr\'e-Oort and Manin-Mumford conjectures. For the latest on the Zilber-Pink conjecture see \cite{pilabook}.

We state here, for the convenience of the reader, a version of the conjecture for curves in Shimura varieties.
\begin{conj}[Zilber-Pink for curves]\label{zpforcurves} Let $\mathcal{X}$ be a Shimura variety, $S_\mathcal{X}$ its set of special subvarieties, and $C\subset \mathcal{X}$ be a curve. Then \begin{center}
		$\atyp(C):=\underset{\underset{\codim(Y)\geq 2}{Y\in S_{\mathcal{X}}}}{\cup} (Y\cap C)$
	\end{center}is a finite set of points assuming that $C$ is not contained in a proper special subvariety of $\mathcal{X}$.
\end{conj}
In the case where $\mathcal{X}$ is some product of copies of $Y(1)$ the conjecture is known to hold conditionally under a large Galois orbits hypothesis by work of P. Habegger and J. Pila, and unconditionally upon enforcing certain conditions on the curve $C$, see \cite{habeggerpila1, habeggerpila2,daworr4} for more. 

Analogous statements can be made for curves in algebraic tori or abelian varieties. In the case of tori the analogue of \Cref{zpforcurves} is known by work of E. Bombieri-P. Habegger-D. Masser-G. Maurin-U. Zannier, see \cite{bmz1,bmz2,bhmz,maurin}.

The aforementioned work of C. Daw and M. Orr follows the general strategy outlined in \cite{dawren}, which can be described as the analogue of the Pila-Zannier method, see \cite{pilazannier}, in the case of Shimura varieties. This strategy reduces the Zilber-Pink conjecture to three hypotheses of arithmetic nature. 

The first of these hypotheses, the "hyperbolic Ax-Lindemann conjecture", has been settled by N. Mok, J. Pila, and J. Tsimerman in \cite{mok2017ax}.  In \cite{daworr, daworr2, daworr3} C. Daw and M. Orr settle cases of the second of these hypotheses, which is a series of conjectures aimed at controlling the ``complexity'' of pre-special subvarieties. Namely, in \cite{daworr,daworr2} the authors establish these hypotheses for all subvarieties of interest in $\mathcal{A}_2$ and in \cite{daworr3} they establish these hypotheses for certain special subvarieties of $\mathcal{A}_g$, namely those of simple PEL type $I$ and $II$. 

The last of the three hypotheses is a so-called ``Large Galois orbits hypothesis'', see for example Conjecture $1.5$ in \cite{daworr3}. In \cite{daworr, daworr2}  C. Daw and M. Orr reduce the Zilber-Pink conjecture for Hodge generic curves in $\mathcal{A}_2$ to the aforementioned hypothesis on Galois orbits, see Theorem $1.2$ of \cite{daworr} and Theorem $1.3$  of \cite{daworr2}. In the more general case of $\mathcal{A}_g$ they obtain analogous results for the intersections of Hodge generic curves in $\mathcal{A}_g$ with simple PEL type $I$ and $II$ special subvarieties, see Theorem $1.3$ of \cite{daworr3}.

They are able to establish the aforementioned hypothesis on Galois orbits for certain curves in $\mathcal{A}_g$. Namely curves defined over $\bar{\Q}$ that intersect the $0$-dimensional stratum of the boundary of the Baily-Borel compactification of $\mathcal{A}_g$. In the case where $g=2$ their results can be summarized as follows: 

\begin{theorem}[\cite{daworr,daworr2}]\label{daworrmain}Let $C$ be an irreducible Hodge generic algebraic curve defined over $\bar{\Q}$ such that the Zariski closure of $C$ in the Baily-Borel compactification of $\mathcal{A}_2$ intersects the $0$-dimensional stratum of the boundary.
	
	Then there are only finitely many points of intersection of $C$ with $E\times{CM}$-curves and quaternionic curves. 
\end{theorem}

Likewise they establish the hypothesis on Galois orbits in the case of such curves in  $\mathcal{A}_g$ for what they call simple PEL type $I$ and $II$ points, see Theorem $1.6$ in \cite{daworr3}.

The key input in establishing the hypothesis on Galois orbits is Andr\'e's G-functions method and the height bounds, as in \Cref{andresbound}, it provides. The condition on the intersection of the curve with the boundary of the Baily-Borel compactification of $\mathcal{A}_g$ also stems from, essentially, the conditions in \Cref{andresbound}.

Following the exposition and results of \cite{daworr3}, using \Cref{application2}, we are able to establish some cases of the Zilber-Pink conjecture unconditionally for certain curves in $\mathcal{A}_g$. The most general case that follows from our height bounds is the following:

\begin{theorem}\label{zpintro}Let $Z\subset \mathcal{A}_g$ be a Hodge generic curve defined over $\bar{\Q}$, where $g\geq 2$. Assume that the Zariski closure of the curve $Z$ in the Baily-Borel compactification of $\mathcal{A}_g$ intersects the stratum of the boundary that is isomorphic to $\mathcal{A}_{g-h}$ with $h\geq 2$.
	
	Then there are only finitely many strongly exceptional points of simple PEL type $I$ and $II$ in $C(\bar{\Q})$.
	\end{theorem}

\begin{rmks}
	1. For the notion of ``strong exceptionality'' in this case see \Cref{stronglyexceptionalag}.\\
	
	2. The aforementioned \Cref{zpintro}, while general, looks perhaps more restrictive than it actually is. We have dedicated the last section of our exposition to presenting some more concrete examples of Zilber-Pink-type statements, see \Cref{section:applicationtozp} that follow, essentially, as corollaries of \Cref{zpintro}.
\end{rmks}

\subsection{Outline of the paper}

 We start, in \Cref{section:survey}, with a short review of the limit mixed Hodge structure associated to a variation of Hodge structures coming from a one-parameter family of abelian varieties that degenerates at some isolated point.

In \Cref{section:abeliangfunctions} we prove two results that we need and concern the behavior of relative periods archimedeanly close to the point $s_0$. In particular, we start in \Cref{section:abeliangfunctionsexistence} by establishing a variant of Theorem $1$ in Ch. IX. \S $4$ of \cite{andre1989g} that guarantees the existence of G-functions among the relative periods archimedeanly close to the point $s_0$. In \Cref{section:independence} we establish the equivalent of Lemma $9.1$ of \cite{papas2022height}, or the Lemma on page $209$ of \cite{andre1989g}, in the case of abelian schemes. This will guarantee that the G-functions that we created are independent of the chosen archimedean place of $K$ considered in their construction in \Cref{section:abeliangfunctionsexistence}.

In \Cref{section:settingnontrivial} we introduce a new setting that is needed to establish the main height bounds. Namely one has to work over a an appropriate cover of the curve $S$. Associated to these covers one gets a family of G-functions coming from the results of the previous section. In this new setting we are able to apply the results of \cite{papas2022height} to define so called polynomial relations among the values of the aforementioned G-functions at points that are archimedeanly close to $s_0$ and satisfy the conditions in \S $8.4.1$ of \cite{papas2022height}. These relations will turn out to be ``non-trivial'', in the notation of \cite{andre1989g}, by results established in $\S7$ of \cite{papas2022height}.

Afterwards, in \Cref{section:globality} we establish ``globality'' for the relations created in \Cref{section:nontrivialrelations} at strongly exceptional points. This is achieved by essentially following the same strategy employed by Andr\'e in his original proof using Gabber's Lemma to establish that the points in question cannot be ``$v$-adically close'' to the degeneration point $s_0$ for any non-archimedean place of their field of definition. Finally, in \Cref{section:everythingtogether} we put everything together and establish \Cref{application2}. We furthermore discuss some interesting examples of strongly exceptional points in \Cref{section:examples}.

In our final section, \Cref{section:caseszp}, we discuss applications of the height bounds we obtain to the Zilber-Pink conjecture. We start by establishing a version of the large Galois orbits hypothesis for simple strongly exceptional points, see \Cref{largego}. This is the crucial ingredient needed, together with the results of \cite{daworr3}, in order to establish \Cref{zpintro} above. We close off our exposition by discussing some specific examples of Zilber-Pink type statements that are significantly more concrete than \Cref{zpintro}. These follow from the more general aforementioned result building on \cite{daworr3} and \cite{tsimermanag}.\\

\textbf{Acknowledgments:} The author thanks Jacob Tsimerman for many helpful discussions. He also thanks Chris Daw for his encouragement, helpful remarks, and for answering some questions about his work with M. Orr. The author finally thanks David Urbanik for some comments on an earlier draft of this paper and Morihiko Saito for making his paper \cite{saitoanf} available. The author was partially supported by Michael Temkin's ERC Consolidator Grant 770922 - BirNonArchGeom.


	\section{General background on limit mixed Hodge structures}

We collect some well known facts about the limit mixed Hodge structure of a $1$-parameter degenerating family of abelian varieties over an algebraically closed subfield $k$ of $\C$. Our main references for this section are \cite{delignehodge3, jumpsmono, saitoanf}.

\subsection{Limit mixed Hodge structures for $1$-parameter families of abelian varieties: A survey}\label{section:survey}

Throughout this text we want to study abelian schemes $f:X\rightarrow S$, where $S$ is the complement of some point $s_0$ in some smooth connected curve $S'$, such that the scheme $f:X\rightarrow S$ degenerates over $s_0$. We are primarily interested in the case where all of the above are defined over some fixed number field $K$, but for certain results that we need it is enough to consider the case where all of the aforementioned data is defined over some algebraically closed field $k$ of characteristic $0$ contained in $\C$. 

In this section for reasons of convenience and notational simplicity we concentrate on the latter case and fix throughout an algebraically closed field $k$ with $char (k)=0$ and an embedding $k\hookrightarrow \C$, as well as $S'$, $s_0$, $S$ and $f:X\rightarrow S$ defined as above and let $g>1$ be the dimensions of the fibers of the morphism $f$.

Associated to this data, after taking the analytification of all objects in question, one obtains several variations of pure polarized Hodge structures of weight  $n$, namely the variations over $S^{an}$ whose underlying local systems is $R^nf^{an}_{*}\Q$. It is classical, see \cite{mumfordabelian}, that the weight $1$  variation ``determines'' the rest of these variations, as they coincide with $\wedge^nR^1f^{an}_{*}\Q$. So from now on we focus on the variation with $n=1$.

Since all of the data is defined over the algebraically closed field $k$, that for all intents and purposes we may think of as the field $\bar{\Q}$ as far as our applications are concerned, we have an associated vector bundle $H^1_{DR}(X/S)$ of $k$-vector spaces over $S$ together with an integrable connection $\nabla$, the Gauss-Manin connection, that is defined over the field $k$ by work of Oda-Katz, see \cite{katzoda}. We also have a vector subbundle $F^1:=e^{*} \Omega^{1}_{X/S}$ of $H^1_{DR}(X/S)$, where $e:S\rightarrow X$ is the zero section of the abelian scheme.

The relative version of the Grothendieck comparison isomorphism between algebraic de Rham and Betti cohomology then takes the form
\begin{equation}\label{eq:relgrothcohomo}
	P_{X/S}:H^1_{DR}(X/S)\otimes_{\mathcal{O}_S}\mathcal{O}_{S^{an}}\xrightarrow{\simeq } R^1(f^{an})_{*}\Q\otimes \mathcal{O}_{S^{an}}.
\end{equation}For this isomorphism we know that $F^1\otimes_{\mathcal{O}_S}\mathcal{O}_{S^{an}}$ defines, under the above isomorphism, the Hodge bundle of the variation on the right hand side of the above isomorphism. For that reason we refer to $F^1$ as the \textbf{Hodge bundle} of the variation as well.

\subsubsection*{Variation of homology}

By Poincar\'e duality we also have a variation of polarized Hodge structures that naturally keeps track of the first homology group of the fibers of $f^{an}$. Namely we have the local system given by \begin{equation}
	(R^{2g-1}f^{an}_{*}\Q)(g)\simeq (R^1f^{an}_{*}\Q)^{\vee},
\end{equation}which is a variation of polarized weight $-1$ pure Hodge structures.

Using this, Grothendieck's comparison isomorphism becomes\begin{equation}\label{eq:relgrothhomology}
		P_{X/S}:H^1_{DR}(X/S)\otimes_{\mathcal{O}_S}\mathcal{O}_{S^{an}}\xrightarrow{\simeq } ((R^{2g-1}(f^{an})_{*}\Q)(g))^{\vee}\otimes_{\Q} \mathcal{O}_{S^{an}}.
\end{equation}

\subsubsection*{Limit mixed Hodge structures}

The above comparison isomorphisms have canonical extensions over $(S')^{an}$ by the notion of a limit mixed Hodge structure at the point $s_0$. We review here the basics of this construction that we will need later on.

Given an abelian scheme $f:X\rightarrow S$ with $g$-dimensional fibers as above we let $f':X'\rightarrow S'$ be its connected N\'eron model over $S'$. We also consider $\bar{f}:\bar{X}\rightarrow S'$ its compactification, by \cite{nagata}, \cite{hironakaI}, and \cite{hironakaII}. By considering base change with an algebraic branched cover of $S'$ we may, for our purposes, a priori assume that the following is true:
\begin{conv}\label{conventinolocalmonodromy}
	The local monodromy around $s_0$ acts unipotently on the fibers of $R^1f^{an}_{*} \Q$ in an analytic neighborhood of $s_0$.
	
\end{conv}

The limit mixed Hodge structures we are interested in have two realizations that are isomorphic via the limit of the isomorphism $P_{X/S}$ above at the point $s_0$. The algebro-geometric one gives us, as in the non-mixed case, an underlying $k$-vector space, while the analytic-geometric one gives us the underlying $\Q$-vector space of the mixed Hodge structure in question.

In more detail, on the one hand, by work of Steenbrink \cite{steenbrink}, we have an algebro-geometric realization of this limit via logarithmic de Rham cohomology. In particular, Steenbrink, see $\S 1$ of \cite{steenbrink}, defines the relative de Rham complex of $\bar{X}$ over $S'$ with logarithmic poles along $\bar{f}^{-1}(s_0)$ that consists of the sheaves $\Omega^{\bullet}_{\bar{X}/S}(\log (\bar{f}^{-1}(s_0)))$ and the hyper-derived direct images $H^p_{DR}(\bar{X}/S'(\log(\bar{f}^{-1}(s_0)))):=\R^pf_{*} \Omega^{\bullet}_{\bar{X}/S}(\log (\bar{f}^{-1}(s_0)))$.

The sheaf $H^1_{DR}(\bar{X}/S'(\log(f^{-1}(s_0)))$ is a coherent sheaf of $\mathcal{O}_{S'}$-modules and is also endowed with the Gauss-Manin connection $\nabla$. Because of \Cref{conventinolocalmonodromy}, this is a connection with nilpotent residue at $s_0$, see \cite{deligneregular}, $II.3.11$. In fact, the pair $(H^1_{DR}(\bar{X}/S'(\log(\bar{f}^{-1}(s_0))), \nabla)$ is Deligne's canonical extension of $H^1_{DR}(X/S)$, see  \cite{deligneregular}, $II.5.2$. We denote this sheaf by $\mathcal{E}$. The Hodge bundle $F^1$ also has an algebro-geometric extension which we will denote by $\mathcal{F}^1$ and coincides with $(e')^{*}\Omega^{1}_{X'/S'}$, where $e'$ is the zero section of the semiabelian scheme $X'\rightarrow S'$.

On the other hand, the limit mixed Hodge structure, or rather its underlying $\Q$-vector space, has a description via the complex of nearby cycles associated to the morphism $\bar{f}^{an}$. Indeed, one has the limit cohomology and homology of $\bar{f}^{an}$ at $s_0$
\begin{center}
	$H^1(X_\infty,\Q):= \mathbb{H}^1(\bar{X}_{s_0}, R\psi_{\bar{f}}(\Q))$, and respectively
\end{center}
\begin{center}
	$H_1(X_\infty,\Q):= \mathbb{H}^{2g-1}(\bar{X}_{s_0}, R\psi_{\bar{f}}(\Q))(g)$.
\end{center}

\subsection{Properties of the limit MHS}\label{section:properties}

With notation as in \Cref{section:survey} we record here some of the main properties of the limit MHS $H^1(X_{\infty},\Q)$ and $H_1(X_{\infty},\Q)$ that we will need later on. The main feature that we will take advantage of is that, in contrast to a general limit of Hodge structures, we have additional structure over the point $s_0$, namely the degeneration of the abelian scheme $f:X\rightarrow S$ at the point $s_0$, as recorded by the fiber over $s_0$ of its N\'eron model over $S'$.\\

With these things in mind, we consider $f:X\rightarrow S$ as above, assume that \Cref{conventinolocalmonodromy} holds and let $X'_0:=X'_{s_0}$ be the fiber over $s_0$ of its connected N\'eron model. In our case $X'_0$ is a semiabelian variety over $k$ and we let $h$ be the dimension of its toric part, namely we assume that the Chevalley decomposition of $X'_0$ is of the form
\begin{equation}\label{eq:chevalley1}
	0\rightarrow T\rightarrow X'_0\rightarrow B\rightarrow 0,
\end{equation}with $T$ an $h$-dimensional torus over $k$ and $B$ a $(g-h)$-dimensional abelian variety over $k$.

\subsubsection*{The limit of cohomology}

As mentioned in \Cref{section:survey}, we have a mixed $\Q$-Hodge structure, which we will denote by $(H^1,W_{\bullet},(F_0^{\bullet})_{\C})$, appearing as the limit of the variation of Hodge structures $R^1f^{an}_{*}\Q$. This mixed Hodge structure is of type \begin{center}
	$\{(0,0),(1,0),(0,1),(1,1)  \}$
\end{center}and we have its weight filtration given by
\begin{center}
	$0=W_{-1}\subset W_0\subset W_1\subset W_2=H^1(X_{\infty},\Q)$, 
\end{center}and its Hodge filtration given by 
\begin{center}
	$0=(F_0^2)_{\C}\subset (F_0^1)_{\C}\subset (F_0^0)_{\C}=H^1_{\C}$
\end{center}

We note that we have the $k$-vector space $\mathcal{F}^1_{s_0}$, the fiber over $s_0$ of the sheaf $\mathcal{F}^1$. Considering the limit of the isomorphism \eqref{eq:relgrothcohomo} over $s_0$, which we denote by $P_0$, we have that $(F_0^1)_{\C}$ is identified with $\mathcal{F}^1_{s_0}\otimes_{k}\C$.

We note that $\mathcal{F}^1_{s_0}$ has a convenient description coming from the fact that our variation comes from a degenerating abelian scheme. In fact we have \begin{equation}\label{eq:f1semiabe}
	\mathcal{F}^1_{s_0}= e_0^{*}\Omega^{1}_{X'_0/k},
\end{equation}where $e_0$ denotes the zero section of the semiabelian scheme $X'_0/k$.

Finally, we note that in view of \eqref{eq:chevalley1} we have a short exact sequence among differential forms\begin{equation}\label{eq:exactdifferent}
	0\rightarrow e^{*}_{0} \Omega^{1}_{B/k}\rightarrow e^{*}_{0} \Omega^{1}_{X'_0/k}\rightarrow e^{*}_{0} \Omega^{1}_{T/k}\rightarrow 0,
\end{equation}where $e_0$ is once again the zero section of these algebraic groups over $k$.

\subsubsection*{The limit of homology}

As mentioned above, we also have a mixed $\Q$-Hodge structure, which we will denote by $(H_1, \widehat{W}_{\bullet}, \widehat{F}^{\bullet}_{0})$, that appears as the limit mixed Hodge structure of the variation $R^{2g-1}f^{an}_{*}\Q(g)$. This will be a mixed Hodge structure of type 
$\{(0,0),(-1,0),(0,-1),(-1,-1)  \}$ and we have its weight filtration given by
\begin{center}
	$0=\widehat{W}_{-3}\subset\widehat{W}_{-2}\subset\widehat{W}_{-1}\subset\widehat{W}_0=H_1(X_{\infty},\Q)$, 
\end{center}and its Hodge filtration given by 
\begin{center}
	$0=(\widehat{F}_0^1)_{\C}\subset (\widehat{F}_0^{0})_{\C}\subset (\widehat{F}_0^{-1})_{\C}=(H_1)_{\C}$
\end{center}

\subsubsection*{Interplay}

The above mixed Hodge structures are connected by the fact that $H^1$ is naturally identified with the dual of $H_1$, see \cite{jumpsmono} Propositions $5.4$ and $6.1$. In particular we have for the two weight filtrations that, under this identification,
\begin{center}
	$W_{j}=(\widehat{W}_{-j-1})^{\perp}$.
\end{center}

For the Hodge filtrations of these mixed Hodge structures we have by duality, see \cite{peterssteenbrink} Example $3.2.2$, a natural identification 
\begin{center}
	$(F^1_0)_\C =(\widehat{F}_{0}^0)^{\perp}_{\C}$.
\end{center}

\subsubsection*{1-motives and the limit homology}

We can extract information about the limit homology from the structure of the semiabelian variety $X'_0$ and its Chevalley decomposition. Our main references are \cite{delignehodge3}, and in particular $\S10.1$, and \cite{jumpsmono}, and in particular $\S 6$.\\

The limit homology appears more naturally as a $\Z$-mixed Hodge structure whose underlying lattice $(H_1)_{\Z}:=H_1(X_{\infty},\Z)$ is such that $H_1(X_{\infty},\Q)=H_1(X_{\infty},\Z)\otimes_{\Z}\Q$. In this guise, we have Deligne's theory of (mixed) $1$-motives. Deligne's results, see \S $10.1$ of \cite{delignehodge3}, describe an equivalence of categories between the category of semiabelian varieties over $\C$ and that of mixed Hodge $1$-motives with negative weights.

Note that, first of all, we have that the analytification of the Chevalley decomposition \eqref{eq:chevalley1} induces a short exact sequence\begin{equation}\label{eq:exacthomology}
	0\rightarrow H_1(T^{an},\Z)\rightarrow H_1((X'_0)^{an},\Z)\rightarrow H_1(B^{an},\Z)\rightarrow 0.
\end{equation}We also have the following, see \cite{delignehodge3} $10.1.3.1$, commutative diagram, with $\Lambda$ a free finite type $\Z$-module,
\begin{equation}\label{eq:commdiagram}
	\begin{tikzcd}
	0 \arrow[r] & H_1((X'_0)^{an},\Z) \arrow[r]                                & \lie ((X'_0)^{an}) \arrow[r, "\exp "']                 & (X'_0)^{an} \arrow[r]                 & 0 \\
	0 \arrow[r] & H_1((X'_0)^{an},\Z) \arrow[ Rightarrow, u, equal] \arrow[r] & (H_1)_{\Z} \arrow[u, "\alpha"'] \arrow[r] & \Lambda \arrow[u, "\beta"'] \arrow[r] & 0.
\end{tikzcd}
\end{equation}

From the above we obtain natural identifications of $\widehat{W}_{-1}(H_1)_{\Z}= H_1((X'_0)^{an},\Z)$ and $\widehat{W}_{-2}(H_1)_{\Z}=H_1(T^{an},\Z)$. In particular we have that $\dim_{\Z} \widehat{W}_{-2}=h$ and $\dim_{\Z}\widehat{W}_{-1}= 2g-h$. One also has that $\widehat{F}_0^0=\ker(\alpha_{\C})$ and $\widehat{F}_0^0\cap \widehat{W}_{-2}=0$.


\section{G-functions as periods}\label{section:abeliangfunctions}

In this section we establish two results, in analogy with results in \cite{andre1989g}, that show that certain relative periods close enough to a degeneration are in fact G-functions. We start with establishing the existence of these G-functions in \Cref{abeliangfunctions}, essentially after we fix an archimedean place. The second results of this section shows that in fact the G-functions in question are in a sense ``independent of this archimedean embedding''.

\subsection{Existence of G-functions}\label{section:abeliangfunctionsexistence}

In this subsection for convenience we adopt slightly different notation from the sections that follow. We let $k=\overline{\Q}$ and $S'$ be a smooth connected curve over $k$ together with a point $s_0\in S'(k)$, and consider $f:X\rightarrow S$ to be a $G$-admissible abelian scheme over $S'$. We also fix an embedding $k\hookrightarrow \C$ and consider a unit disk $\Delta\subset S'^{an}_{\C}$ centered at the point $s_0$ so that for the punctured disc $\Delta^{*}$ we have that $\Delta ^{*}\subset S$. 

We also let $x$ be a local parameter of $S'$ at the point $s_0$ and let $X'$, $T$, $B$, etc. be as in \Cref{section:notation}.
As in \Cref{section:survey} we assume that \Cref{conventinolocalmonodromy} holds so that local monodromy around the point $s_0$ is unipotent.

Our aim here is to establish an analogue of Theorem $1$ on pages $184$-$185$ of \cite{andre1989g}.
\begin{theorem}\label{abeliangfunctions}
	For the above objects there exists a basis $\omega_1,\ldots,\omega_{2g}$ of $H^1_{DR}(X/S)$ over some dense open subset $U$ of $S$, such that the Taylor expansion in $x$ of the relative periods $\frac{1}{2\pi i} \int_{\gamma_j}^{} \omega_i$ are G-functions for $1\leq j\leq h$, where $\gamma_j$ constitute a local frame of the sheaf $M_{0}\R_1 (f^{an}_{\C})_{*}(\Q)|_{V}$, where $V\subset \Delta^{*}$ is some open analytic subset.
\end{theorem}

\begin{rmk}
Here $M_{0}\R_1 (f^{an}_{\C})_{*}(\Q)$ is the image of the logarithm of the local monodromy around $s_0$ acting on the sheaf $R^{2g-1}(f^{an})_{*}\Q|_{\Delta^{*}}$. We borrow this notation from Andr\'e, see page $185$ of \cite{andre1989g}. 
\end{rmk}

\subsection{Proof of \Cref{abeliangfunctions}}


	
	
		
		
	
	
	We follow closely the proof of Theorem 1, Ch. X, \S 4 of \cite{andre1989g} which we are trying to generalize.
	\begin{proof}
	
	We consider the coherent $\mathcal{O}_{S}$-module with connection $(H^1_{DR}(X/S),\nabla)$. By Deligne's work on differential equations with regular singular points \cite{deligneregular} and Griffiths' regularity theorem, it is known that there exists an extension, i.e. a locally free coherent $\mathcal{O}_{S'}$-module, $\mathcal{E}$ of this sheaf. As we mentioned in \Cref{section:survey} this has an algebro-geometric description as the relative logarithmic de Rham cohomology sheaf $H^1_{DR}(\bar{X}/S'(\log(\bar{f}^{-1}(s_0))))$, where $\bar{f}:\bar{X}\rightarrow S'$ is, as in \Cref{section:survey}, a compactification of the connected N\'eron model.
	
	It is also known that $\mathcal{F}^1:=e^{*}\Omega^{1}_{X'/S'}$ is a subbundle of $\mathcal{E}$ and that $(e^{*}\Omega^{1}_{X'/S'})|_{S}$ is nothing other than the Hodge bundle $F^1$. In other words the analytification of the complexification of $\mathcal{F}^{1}$ defines the Hodge filtration fiber-wise on $H^1(X^{an}_{s},\Q)$ for all $s\in S^{an}$.
	
	We start by establishing the following claim, in parallel to the proof of Andr\'e's result in loc.cit.:
	
	\begin{claim}\label{claim1}Let $\{\gamma_i:1\leq i\leq h\}$ be a $\Q$-basis of the vector space $H_1(T^{an},\Q)$. There exists a basis $\{w_i:1\leq j\leq 2g\}$ of the $k$-vector space $(\mathcal{E})_{s_0}$ such that $\frac{1}{2\pi i} \int_{\gamma_j}^{} w_i\in k$ for all $1\leq i\leq 2g$ and $1\leq j\leq h$. Equivalently, $\frac{1}{2\pi i} \int_{\gamma}^{} \omega\in k$ for all $\omega\in (\mathcal{E})_{s_0}$ and all $\gamma\in H_1(T^{an},\Q)$.
	\end{claim}

\textbf{Notation: }We introduce a slight change in our notation for reasons of convenience, we will denote the fibers of any object $(*)$ at the point  $s_0$ by writing $(*)_0$ instead of $(*)_{s_0}$.\\

\begin{proof}[Proof of \Cref{claim1}] We have a natural isomorphism, given by the limit of the comparison isomorphism $P_{X/S}$ at $s_0$, 
\begin{center}
		$P_0: \mathcal{E}_0\otimes_k \C\xrightarrow{\simeq} (H_1(X_{\infty}, \Q)\otimes_{\Q} \C)^{\vee}$,
\end{center}where $H_1(X_{\infty},\Q)$ is the limit mixed Hodge structure of the variation of Hodge structures whose underlying local system is $R^{2g-1}(f^{an})_{*}(\Q)(g)$. As noted, we may and do identify the fibers of this variation with the first homology group of the respective fiber. 
	
	Consider the fiber $X'_0$, which is a semiabelian variety with Chevalley decomposition given over $k$ by 
	\begin{equation}\label{eq:chevalley}	
		0\rightarrow T \rightarrow X'_{0}\rightarrow B\rightarrow 0,
	\end{equation}with $T$ an $h$-dimensional $k$-torus and $B$ a $(g-h)$-dimensional abelian variety over $k$. 

This decomposition induces an exact sequence between the first homology groups of the analytifications of the algebraic groups in question
\begin{equation}\label{eq:exhom}
	0\rightarrow H_1(T^{an},\Q)\rightarrow H_1((X'_{0})^{an},\Q)\rightarrow H_1(B^{an},\Q)\rightarrow 0.
\end{equation}We also have a short exact sequence among the pullbacks of differentials
\begin{equation}\label{eq:exdif}
	0\rightarrow e_0^{*} \Omega^{1}_{B/k} \rightarrow e_0^{*} \Omega^{1}_{X'_0/k}\rightarrow e_0^{*} \Omega^{1}_{T/k}\rightarrow 0,
\end{equation}via the zero section $e_0$ of each algebraic group.

With this notation, we consider a $k$-basis $\{\beta_i:1\leq i\leq h\}$ of $e^{*}_{0}\Omega^1_{T/k}$ and for each $i$ let $w_i\in e^{*}_{0} \Omega^{1}_{X'_0/k}$ be such that $w_i\mapsto \beta_i$ under the rightmost arrow in \eqref{eq:exdif}. We extend the set $\{w_i:1\leq i\leq h\}$ to a basis $\{ w_i:1\leq i\leq g\}$ of $ e^{*}_{0}\Omega^{1}_{X'_0/k}$ by choosing $w_i\in e^{*}_{0}  \Omega^{1}_{B/k}$ for $i\geq h+1$ that form a basis of $ e^{*}_{0}\Omega^{1}_{B/k}$.
	
By the unipotence of the local monodromy around $s_0$, we know that the weight filtration of the limit MHS $H_1(X_\infty,\Q)$ will have the form \begin{equation}\label{eq:homweight}
	0=\widehat{W}_{-3}\subset \widehat{W}_{-2}\subset \widehat{W}_{-1}\subset \widehat{W}_0=H_1(X_{\infty},\Q)
\end{equation}and that we have $\widehat{W}_{-1}=H_1((X'_0)^{an},\Q)$ and $\widehat{W}_{-2}=H_1(T^{an},\Q)$. 
We let $\{\gamma_j:1\leq j\leq h\}$ be a $\Q$-basis of $\widehat{W}_{-2}$.

The non-degeneracy of the pairing $H^0(T,\Omega^1)\times H_1(T^{an},\Q)\rightarrow k$ induced by the residue, guarantees that the $h\times h$ matrix $(\frac{1}{2\pi i} \int_{\gamma_j}^{} \omega_i)_{1\leq i,j\leq h}$ is in $\GL_h(k)$. 

For the $w_i$ with $h+1\leq i\leq g$ we claim that in this case $\frac{1}{2\pi i} \int_{\gamma_j}^{} \omega_i=0$ for $1\leq j\leq h$. 

This can be seen as follows:\\

From the mixed Hodge structure $((H_1)_{\Z},\widehat{W}_{\bullet},\widehat{F}_{0})$, with notation as in \ref{section:properties}, we get a short exact sequence of Jacobians, in the notation of \cite{saitoanf} $1.1.5$, see also Lemma $10.1.3.3$ of \cite{delignehodge3},
\begin{equation}\label{eq:sesjacobians}
	0\rightarrow J( Gr^{\widehat{W}}_{-2}H_1)\rightarrow J(H_1) \rightarrow J(Gr^{\widehat{W}}_{-1}H_1)\rightarrow 0
\end{equation}that is nothing but the analytification of the Chevalley decomposition of $X'_{s_0}$.

This gives a short exact sequence among the respective Lie algebras of these analytic groups
\begin{equation}\label{eq:sesliealg}
	0\rightarrow \lie T^{an}\rightarrow \lie (X_0')^{an}\rightarrow \lie B^{an}\rightarrow 0.\end{equation}

This is naturally the dual of the exact sequence of the analytification of the complexification of \eqref{eq:exdif}. From this duality we get that if $w\in e^{*}_{0}\Omega^{1}_{B^{an}/\C}$, in particular if $w\in e^{*}_{0}\Omega^{1}_{B/k}$, and $\gamma\in H_1(T^{an},\Q)\subset \lie T^{an}$ we have $\frac{1}{2\pi i}\int_{\gamma}w=0$.\\

All that is left, therefore, is to extend the basis $\{w_i :1\leq i\leq g\}$ of $\mathcal{F}^1_0$ to one of $\mathcal{E}_0$ and show that these new elements of the basis satisfy $\frac{1}{2\pi i}\int_{\gamma_j}w\in k$.

For that we follow the notation in \Cref{section:survey} with $\bar{f}:\bar{X}\rightarrow S'$ a compactification of the connected N\'eron model over $S'$. In this setting, and under \Cref{conventinolocalmonodromy}, the fiber $\bar{X}_0:=\bar{f}^{-1}(s_0)$ is a reduced divisor with normal crossings. We have thus the description 
$\mathcal{E}=H^1_{DR}(\bar{X}/S'(\log(\bar{X}_0)))$ and by Proposition $2.16$ of \cite{steenbrink} we have that the fiber $\mathcal{E}_0$ is canonically isomorphic to $\mathbb{H}^1(\bar{X}_0,\Omega^{\bullet}_{\bar{X}/S'}(\log(\bar{X}_0))\otimes_{\mathcal{O}_{\bar{X}}}\mathcal{O}_{\bar{X}_0})$.  

\begin{claim}\label{claim2}We have a natural map $\pi:\mathcal{E}_0 \rightarrow H^1_{DR}(T/k)$ whose complexification coincides under the comparison isomorphism $P_0$ with the map \begin{center}
	$H^1(X_{\infty},\C)\rightarrow Gr^{W}_{2}   H^1 =   H^1(T^{an},\C)\rightarrow 0$.
\end{center}
\end{claim}

Assuming \Cref{claim2}, we are done. Indeed, letting $P_0(w)=\Sum{j=1}{h} a_j  \gamma_j^{*} +W_{1,\C}$, where $\gamma_j^{*}$ are the dual elements of the $\gamma_j$, we get that \begin{center}
	$\frac{1}{2\pi i}a_j=\frac{1}{2\pi i}    \int_{\gamma_j}^{} \pi(w) $
\end{center}which is in $k$, as it coincides with the residue of $\pi(w)$ at $0$, similarly to our previous argument.
	
\end{proof}

\begin{proof}[Proof of \Cref{claim2}]We write $i:X'\hookrightarrow \bar{X}$ for the open immersion of the connected N\'eron model $X'$ of $X$ over $S'$ in its compactification $\bar{X}$ that we chose earlier. We note that $f':X'\rightarrow S'$ is smooth by definition of the N\'eron model.
	
First of all by pullback, via the open immersion $i_0:X'_0\hookrightarrow \bar{X}_0$, we get a canonical map \begin{equation}\label{eq:pullbackmap}
	\mathbb{H}^1(\bar{X}_0,\Omega^{\bullet}_{\bar{X}/S'}(\log \bar{X}_0)) \rightarrow 	\mathbb{H}^1(X_0',\Omega^{\bullet}_{X'/S'}(\log X'_0)),
\end{equation}of $k$-vector spaces. 

On the other hand the smoothness of $f'$ implies that\begin{equation}\label{eq:lognotneeded}
		\mathbb{H}^1(X_0',\Omega^{\bullet}_{X'/S'}(\log X'_0))= \mathbb{H}^1(X'_{0},\Omega^{\bullet}_{X'_0/k})=H^1_{DR}(X'_0/k).
\end{equation}Furthermore, from the inclusion $T\hookrightarrow X'_0$ we get a natural map, again via pullback, \begin{equation}\label{eq:pullbacktorusderham}
H^1_{DR}(X'_0/k)\rightarrow H^1_{DR}(T/k).
\end{equation}

Combining all of the above we end up with a canonical map 
\begin{equation}\label{eq:algderhamatlimit}
\pi:\mathcal{E}_0\rightarrow H^1_{DR}(T/k)\rightarrow 0.
\end{equation}

The complexification of \eqref{eq:algderhamatlimit} coincides under the isomorphism $P_0$ with the map\begin{center}
	 $H^1(X_{\infty},\C)\rightarrow Gr^W_2 H^1=H^1(T^{an},\C)\rightarrow 0$.
\end{center}Indeed, by canonicity, the composition of \eqref{eq:pullbackmap} and \eqref{eq:lognotneeded} coincides, under $P_0$ with the dual of the map $H_1((X_0')^{an},\Z)\otimes_{\Z}\C\rightarrow (H_1)_{\Z}\otimes_{\Z}\C$ appearing at the bottom row of \eqref{eq:commdiagram}. On the other hand, the complexification of the map \eqref{eq:algderhamatlimit} coincides with the dual of the inclusion 
\begin{center}
	$H_1(T^{an},\C)\hookrightarrow H_1((X'_0)^{an},\C)$.
\end{center}
\end{proof}

\begin{rmks}1. We note that the above can be rephrased with respect to the weight filtration on $H^1(X_\infty,\Z)$. Indeed, by the exposition in \Cref{section:survey}, we have that the the complexification of the maps $\mathcal{E}_0\rightarrow H^1_{DR}(X'_0/k)$ and $\mathcal{E}_0\rightarrow H^1_{DR}(T/k)$ in the above proof correspond under $P_0$ with the projections \begin{center}
		$H^1(X_{\infty},\C)\rightarrow H^1(X_{\infty},\C)/W_{0,\C}$, and 
	\end{center}
\begin{center}
	$H^1(X_{\infty},\C)\rightarrow H^1(X_{\infty},\C)/W_{1,\C}$ respectively,
\end{center}where $W_{i,\C}:=W_i\otimes_{\Q}\C$.\\

2. We note that all the analytifications in this proof are happening with respect to our fixed embedding $k\hookrightarrow \C$. We suppress mention of this embedding for notational simplicity but we return to this issue in \Cref{abelianchangeofplace}. 
\end{rmks}

We now return to the proof of \Cref{abeliangfunctions}
Consider the basis $\{w_i:1\leq i\leq 2g\}$ of $(\mathcal{E})_{0}$ of the claim.  Nakayama's Lemma allows us then to extend this basis to a basis of sections $\omega_1,\ldots,\omega_{2g}$ of $\mathcal{E}$ over some Zariski neighborhood $U_0:=U\cup \{s_0\}$ of $s_0$, here $U:=U_{0}\backslash\{s_0\}$ is a dense Zariski open subset $S$. In other words, we have that $\omega_i(0)=w_i$ for all $i$.\\

The $\gamma_i$ are monodromy invariant, so they extend to sections of the local system $R^{2g-1}(f^{an})_{*}\Q(g)$ over the punctured disc $\Delta^{*}$. For any $V\subset U^{an}\cap \Delta^{*}$ simply connected open, in the analytic topology of $S'^{an}$, we may extend the set $\{\gamma_i:1\leq i\leq h\}$ to a frame of $R^{2g-1}(f^{an})_{*}\Q(g)|_{V}$. 

With these choices we get the relative period matrix $P:=(\frac{1}{2\pi i} \int_{\gamma_j}^{}\omega_i)$ which will be such that the matrix $P_h$, that consists of the first $h$ columns of the relative period matrix $P$, will be such that it extends at $0$ and furthermore
\[ \vphantom{
	\begin{matrix}
		\overbrace{XYZ}^{\mbox{$R$}}\\ \\ \\ \\ \\ \\
		\underbrace{pqr}_{\mbox{$S$}}
\end{matrix}}%
P_h(0)=\begin{bmatrix}
	A\\
	0\\
	\vdots\\
	0\\
	B
\end{bmatrix}%
\begin{matrix}
	\coolrightbrace{A}{h\text{ rows}}\\
	\coolrightbrace{0\\
		\vdots\\
		0 }{g -h \text{ rows}}\\
	\coolrightbrace{B}{g\text{ rows}}
\end{matrix}\]
with $A\in GL_h(k)$ and $B\in M_g(k)$.\\

Let $k\{x\}$ be the henselization of the local ring $k[x]_{(x)}$ and let $\partial=x\frac{d}{dx}$.
By the theory of the Gauss-Manin connection it is well known that $\nabla_{\partial }$ is represented by a matrix in $M_{2g}(k\{x\})$, with respect to the basis $\omega_i$ of $\mathcal{E}$ that we chose.

We also know that\footnote{See \cite{andre1989g} Ch. $IX$ and in particular the discussion in $\S$ $1.2$ of the aforementioned chapter.} $P_h$ satisfies the differential system \begin{equation}\label{eq:differentialequation}
	\partial X=G\cdot X,
\end{equation}as does the full relative period matrix $P$ and all of its columns.\\

By \Cref{conventinolocalmonodromy} we get that $G(0)$, which coincides with the residue of the connection at $0$, is nilpotent. This allows us to establish the following claim.\\

\begin{claim}
	There exists a matrix $D\in M_{2g\times h}(k)$ such that $P_h$ can be written as \begin{equation}\label{eq:claimandreproof}
		P_h=P_G\cdot D,
	\end{equation}where $P_G:=Y_G \exp(G(0)\log x )$ and $Y_G$ is the normalized uniform part of the solution of \eqref{eq:differentialequation}. 
\end{claim} 

\begin{rmk}We introduce a bit of convenient notation. For any $2g\times 2g$ matrix $M\in M_{2g}(F)$, where $F$ is any field, we will denote by $M_h$ the $2g\times h$ matrix with entries in the field $F$ that consists of the first $h$ columns of the matrix $M$.
\end{rmk}

\begin{proof}[Proof of Claim]We quickly review some notions from Ch. III of loc.cit. The system \eqref{eq:differentialequation} is equivalent to the two systems 
	\begin{equation}\label{eq:equivalentsystems}
		\left.
		\begin{aligned}
			\partial (Z^{-1}X)=CZ^{-1}X  \\
			\partial Z= G\cdot Z-Z\cdot C
		\end{aligned}
		\right\}.
	\end{equation}
	
	In other words $X= Z\exp(C \log x )$, from the first of the two systems. The fact that $G(0)$ is nilpotent, and will hence have weakly prepared eigenvalues\footnote{See \cite{kedlaya} for this definition.}, allows us to choose $C=G(0)$. The normalized uniform part $Y_G$ of $X$ is then the unique solution in $\GL_{2g}(k[[x]])$ of the second system in \eqref{eq:equivalentsystems} that satisfies $Y_G(0)=I_{2g}$.
	
	The existence of the normalized uniform part and the fact that $Y_G\in\GL_{2g}( k[[x]])$ follow from the discussion in Ch. III, \S $1.4$ of loc. cit.\\
	
	Let us set $C:=G(0)$ from now on and define $Z:= P\exp(-C\log x)$. We then have that $Z$ satisfies the second equation in \eqref{eq:equivalentsystems} as does $Y_G$. From the discussion in Ch. III, \S $1.2$ of loc. cit. we get that, since $C$ has weakly prepared eigenvalues, there exists a matrix $D\in M_{2g}(\C)$ such that \begin{center}
		$Z=Y_G \cdot D$. 
	\end{center}
	
	We thus have that $P= Y_G \cdot D \cdot \exp(C \log x)$. We will then have that $P_h= Y_G \cdot D \cdot(\exp(C \log x))_h$.\\
	
	We know that $C=G(0)$ is the residue of the Gauss-Manin connection at $0$ and since the logarithm of the local monodromy has rank $h$ so will $G(0)$, see \cite{deligneregular} II.$3.11$.
	
	Choosing the basis of sections $\omega_i$ appropriately we may arrange so that $G(0)_h$ will be the zero $2g\times h$ matrix. Indeed, computing the logarithm of the local monodromy around the point $s_0$ using the analytic frame $\gamma_j$ one sees that, by choice of this frame and the interplay of this logarithm with the weight filtration, that the logarithm of the local monodromy will have this property.
	
	Hence it suffices to note that we are working over an algebraically closed field $k$ and thus can choose the $\omega_i$ appropriately since $G(0)$ coincides with the logarithm of the monodromy acting on the limit MHS. We will then have that $(\exp(C \log x))_h=I_h +C_h=I_h$ so that, putting everything together, we get that  $P_h= Y_G \cdot D_h$. Since $Y_G(0)=I_{2g}$ we get that $D_h=P_h(0)$ which by our earlier discussion has in fact entries in $k$.	
\end{proof}

The previous claim guarantees that $P_h\in M_{2g\times h}(k[[x]])$. As in the proof of Ch.IX, \S $4.4$ in loc. cit., we get that these entries are solutions to geometric differential equations, as these are defined in Ch.II of loc.cit.. By the above discussion they can also be written as power series with coefficients in $k=\overline{\Q}$. These facts allow us to establish that the entries in question are G-functions by appealing to the Theorem in the appendix of Ch. V in \cite{andre1989g}.
\end{proof}

\begin{rmk}We note that the aforementioned result already appears in the literature, see for example the discussion in $III$, $\S 2$ of \cite{andrepadic}. We chose to include the proof of this result because of its more elementary nature and because it helps in establishing \Cref{abelianchangeofplace} in a very straightforward manner.
\end{rmk}
\subsection{Independence from archimedean embedding}\label{section:independence}

We now return to the setting we will be most interested in. Namely we assume that everything, meaning the curve $S'$ and the semiabelian scheme $f':X'\rightarrow S'$, in the previous subsection is defined over a fixed number field $K$.

 If we were to fix an archimedean place $v_0\in \Sigma_{K,\infty}$, thus getting a corresponding embedding $\iota_0:K\hookrightarrow \C$, \Cref{abeliangfunctions} provides us with power series\footnote{We note that technically here these will have coefficients in some finite extension of $K$ and we might have to tensor our geometric objects with this finite extension. Thus we assume without loss of generality that these coefficients are in fact in $K$ by repeating the process after base changing with said finite extension of $K$ from the start.}  $y_{i,j}(x)\in K[[x]]$, where $1\leq i\leq 2g$ and $1\leq j\leq h$, such that the power series $i_0(y_{i,j})(x)\in \C_v[[x]]$ coincide with the elements of the first $h$-columns of the relative period matrix of $f:X\rightarrow S$, over any $V\subset \Delta_{v_0,R_{v_0}}^{*}$. In other words we have that $y_{i,j}(x(P))$ will be the value of the corresponding element of the stalk at $P\in \Delta_{v_0,R_{v_0}}^{*}$ of the relative period matrix. Here, as usual, $\Delta_{v_0,R}$ is a small enough analytic disk on $(S'\times_{i_0}\C)^{an}$ centered at $s_0$, in particular it is enough to consider $R_{v_0}=\min\{1, R_{v_0}(\iota_0(y_{i,j})(x))\}$ where $R_{v_0}(\iota_0(y_{i,j})(x))\}$ denotes the minimum of the radii of convergence of the power series $\iota_0(y_{i,j})(x)$.

Following the exposition in \cite{andre1989g}, Ch. $X$, $\S 3$, we need to make sure that the G-functions that appear as these ``monodromy invariant'' relative periods, make sense as relative periods analytically locally near $s_0$ irrespective of the chosen archimedean place $v$.

Let us introduce some notation. We write $Y=(y_{i,j}(x))\in M_{2g\times h}(K[[x]])$ for this matrix of G-functions and for any embedding $\iota:K\hookrightarrow \C$ we write $Y_\iota$ for the matrix in $M_{2g\times h}(\C)$ obtained from $Y$ under the embedding $\iota$. We also write $f_{\iota,\C}$ for the base change of the morphism $f:X\rightarrow S$ via the morphism $\iota^{*}:\spec(\C)\rightarrow \spec(K)$ induced by the embedding $\iota$.

For example, note that by construction, we know that $Y_{\iota_0}$ appears as the first $h$ columns of the relative period matrix coming from the isomorphism \begin{center}
	$H^1_{DR} (X/S)\otimes_{\mathcal{O}_S}\mathcal{O}_{S^{an}}|_V\xrightarrow{\simeq} R^1(f_{\iota_0,\C})^{an}_{*}\Q\otimes \mathcal{O}_{S^{an}}|_V            $
\end{center}where $V\subset  \Delta_{v_0,R_{v_0}}^{*}$ is some open analytic subset of a small punctured unit disk centered at $s_0$ as above.

We then have the following analogue of the Lemma on page $209$ of \cite{andre1989g}:

\begin{lemma}\label{abelianchangeofplace}
	Let $\iota$ be a complex embedding of $K$. Then $Y_{\iota}$ is again the matrix that consists of the first $h$ columns of the period matrix attached to the same basis of $\mathcal{E}|_U$ and to some local frame of $(R^{2g-1}(f_{\iota,\C})^{an}_{*}\Q)(g)$ over any simply connected analytic open $V\subset \Delta^{*}_{v,R_{v}}$, where $v$ is the place of $K$ that corresponds to $\iota$ and $\Delta^{*}_{v,R_{v}}$ is defined as above.
\end{lemma}

\begin{rmk}
	For convenience we introduce a last bit of notation. For any $K$-scheme, say $Z$, and for any sheaf, say $\mathcal{F}$, over such a scheme, we write $Z_{\iota}$, respectively $\mathcal{F}_{\iota}$, for their base change via $\iota^{*}:\spec(\C)\rightarrow \spec(K)$. We denote sections of such sheaves by $\omega^{\iota}$, mainly in an attempt to keep tabs on the embedding $\iota$.
\end{rmk}

\begin{proof}We have already chosen $\{\omega_i\}^{2g}_{i=1}$ a basis of $\mathcal{E}_U$ and write $\omega^{\iota}_i$ for their image in $\mathcal{E}_\iota:=\mathcal{E}\otimes_{\iota} \C$, which is naturally identified with $H^1_{DR} (\bar{X}_{\iota}/S'_{\iota}(\log \bar{X}_0))$.
	
	Once again, as in our earlier discussion, we choose $\gamma^{\iota}_j$, $1\leq j\leq h$ sections of $(R^{2g-1}(f_{\iota,\C})^{an}_{*}\Q)(g)|_V$, with $V\subset U^{an}_{\iota}$ as above, such that $(\gamma^{\iota}_j)_0$ are a basis of the $\Q$-vector space $\hat{W}_{-2}H_1((X_\iota)_\infty,\Q)$. In particular, as in the proof of \Cref{abeliangfunctions} each $(\gamma^{\iota}_j)_0$ is identified with the $j$-th component of the torus $T_\iota^{an}\simeq (\C^{\times})^h$. We then extend the set of $\gamma^{\iota}_j$ to a local frame $\{\gamma^{\iota}_j\}^{2g}_{j=1}$ on the simply connected set $V$.
	
	These choices, as we have seen give a relative period matrix which we denote by $P_{\iota}=(z_{i,j})$. In contrast, we write simply $P$ for the original period matrix obtained via the original embedding $\iota_0$.
	
	By the proof of \Cref{abeliangfunctions} we know that the $z_{i,j}$ with $1\leq j\leq h$ at $0$ are residues of the basis $\{(\omega^{\iota}_i)_0\}$ of the fiber $(\mathcal{E}_{\iota})_0$. In particular we have that $z_{i,j}(0)=\frac{1}{2\pi i} \int_{(\gamma^{\iota}_j)_0}(\omega^{\iota}_i)_0$ which is equal to $\iota(y_{i,j})(0)$ by conjugation.
	
	On the other hand we know that the matrix $(P_{\iota})_h$ satisfies the differential system \eqref{eq:differentialequation}, which has coefficients in $K[[x]]$. The proof of the second claim in the proof of \Cref{abeliangfunctions} shows that $(P_{\iota})_h=\iota(Y_G) (P_{\iota})(0)$ which is nothing but $Y_{\iota}$, again by the proof of the aforementioned claim.
\end{proof}
	
		\section{Non-trivial relations}\label{section:nontrivialrelations}

In order to apply \Cref{abeliangfunctions} in practice we need to work in an appropriately chosen cover of our original curve $S'$. The need to work with a finite cover arises from the fact that the uniformizing parameter $x$ at $s_0$ may not have only one zero, or may not have simple zeroes! 

Assuming that the parameter $x$ had only one simple zero, namely $s_0$, then for every place $v\in \Sigma_K$ we would get that the set of points $P\in S'(\C_v)$ with $|x(P)|_v<R$ naturally corresponds with the set of points of an analytic disk $\Delta_{v,R} \subset S'^{an}$ centered at $s_0$. In other words, in this case ``$v$-adic proximity'' of $P$ to the point $s_0$, meaning that $P$ is in some disk as above centered at $s_0$, is equivalent to ``$v$-adic proximity'' of $x(P)$ to $0$. Note that here by $S'^{an}$ we mean either the classic analytification of $S'$ with respect to $v$ when $v$ is archimedean, or the adic space associated to $S'\times_K K_v$.

In order to address this issue we change the setting that appears in the results discussed in the introduction, see e.g. the discussion before \Cref{application2}. In fact we will first prove the main theorems, e.g. the height bound for exceptional points and the Large Galois orbits hypothesis, for a certain good cover of the original curve that appears in \Cref{application2}.\\

The main goal of this section is to first of all define this new setting, that we eluded to above, and fix the notation that will be used for the rest of the exposition. The starting point of what is to follow is the definition of a natural collection of G-functions associated to an abelian scheme over the cover in question, and most importantly a description of so called ``non-trivial'' relations among the values of these G-functions at what we call ``strongly exceptional points'' of the abelian scheme.

\subsection{Setting}\label{section:settingnontrivial}

We begin by fixing the notation that we will use in the sequel. We follow closely the exposition in Chapter $X$, $\S3$ of \cite{andre1989g}. We have also found extremely helpful the exposition in $\S 6$ of \cite{davidg} and $\S 5$ of \cite{daworr4}.\\

Consider $K$ a fixed number field, $S'$ a smooth geometrically connected curve defined over $K$, $s_0\in S'(K)$ a closed point, and $S:=S'\backslash\{ s_0\}$. 

As usual, we consider $f:X\rightarrow S$ an abelian scheme with $g$-dimensional fibers, where $g\geq 2$, let $X'$ be its connected N\'eron model over $S'$, and assume that the fiber $X'_0$ over $s_0$ is a semiabelian variety whose Chevalley decomposition is given by the short exact sequence 
\begin{equation}\label{eq:chevalleysettingsection4}
	0\rightarrow T\rightarrow X'_0\rightarrow B\rightarrow 0
\end{equation}with $T$ an $h$-dimensional torus over $K$, where $h\geq 2$, and $B$ a $(g-h)$-dimensional abelian variety over $K$.

\begin{conv}\label{conventionsemistable} By possibly enlarging the base field $K$, from the start, we may and do assume from now on that the torus $T$ is split over $K$ and that the abelian variety $B$ has everywhere semistable reduction.
	
\end{conv}

Finally, we assume that the image of the morphism $S\rightarrow \mathcal{A}_g$ induced from $f$ is Hodge generic in $\mathcal{A}_g$.

\subsubsection{Good Covers}\label{section:goodcovers}

Now let us consider $\bar{S}'$ the completion of $S'$, so that $\bar{S}'$ is a smooth projective curve. Then Lemma $5.1$ of \cite{daworr4} shows that there exists, after possibly extending the constant field $K$ to some finite extension of our original $K$, a triple $(C_4,x,c)$, where \begin{itemize}
	\item $C_4$ is a smooth projective curve,
	
	\item $x\in K(C_4)$ is a non-constant rational function on $C_0$, and 
	
	\item $c:C_4\rightarrow \bar{S}'$ is a non-constant morphism,

\end{itemize}such that the following are satisfied:\begin{enumerate}
	\item the zeros $\{\xi_1,\ldots,\xi_l\}\subset C_4(\bar{K})$ of $x$ are all simple,
	
	\item for all $i$ we have that $c(\xi_i)=s_0$, and 
	
	\item the induced map $x:C_0\rightarrow \mathbb{P}^1$ is a finite Galois covering.
\end{enumerate}

To streamline our exposition a bit we introduce the following: 
\begin{defn}\label{goodcoversdef}Let $(C_4,x,c)$ be a triple associated to $S'$ as above. Then we call $C_4$ a \textbf{good cover} of the curve $S'$. 
\end{defn}

Now let us set $C':=c^{-1}(S')$ and $C:=c^{-1}(S)$. Pulling back the semiabelian scheme $f':X'\rightarrow S'$ we get a semiabelian scheme $f'_C:X'_C:=X'\times_{S'} C'\rightarrow C'$ and an abelian scheme $f_C:X_C:=X\times_S C\rightarrow S$.

By definition of $f_C$ and $f'_C$, together with the assumptions on $f$ and $f'$ in the previous subsection, we then have that the for the semiabelian scheme $X'_C\rightarrow C'$ and for all roots $\xi_i$ of $x$ there is an isomorphism $X'_{\xi_i}\simeq X'_{s_0}$. In particular, the Chevalley decompositions of the $X'_{\xi_i}$ all coincide with that of $X'_{s_0}$. Furthermore, the map $C\rightarrow \mathcal{A}_g$ induced from $f_C$ has image that is a Hodge generic curve in $\mathcal{A}_g$.

\subsubsection{Associated families of G-functions}	 

Fix $f':X'\rightarrow S'$, $f'_C:X'_C\rightarrow C'$, $S$, and $C$ coming from a good cover $(C_4,x,c)$ of $S'$ as above. Then associated to these we have a natural family of G-functions. 

Indeed from \Cref{abeliangfunctions} we get for each root $\xi_i$, $1\leq i\leq l$ of the uniformizing parameter $x$ a family $\mathcal{Y}_{i}\subset \bar{K}[[x]]$ of G-functions. We note that upon tensoring everything from the beginning with some finite extension $L/K$ we may and do assume from now on that the G-functions in the various families are in fact power series with coefficients in the field $K$.

\begin{defn}\label{associatedgfunctions}Given $(C_4,x,c)$ a good cover as above we define its \textbf{associated family of G-functions} to be the set \begin{center}
		$\mathcal{Y}:=\mathcal{Y}_1\sqcup \ldots \sqcup \mathcal{Y}_l$,
	\end{center}where $\mathcal{Y}_i$ are the above sets of G-functions associated to each root of the uniformizing parameter $x$.
	
\end{defn}

\subsubsection*{An integral model}

Let $f':X'\rightarrow S'$ be as above, $C_4$ a good cover of $S'$, $C'$ and $C$ the (possibly affine) curves described above, and $f'_C:X'_C\rightarrow C'$ the associated semiabelian scheme. Following the discussion in \cite{andre1989g}, Ch. $X$, $\S 3.1$, we fix a regular projective model $\tilde{C}$ of the curve $C'$ over the ring of integers $\mathcal{O}_K$. For $s\in C'(L)$, with $L/K$ finite, we write $\tilde{s}$ for the section of the arithmetic pencil $\tilde{C}\times_{\spec(\mathcal{O}_K)}\spec(\mathcal{O}_L)$ induced from the point $s$.

For the G-functions that comprise the family $\mathcal{Y}$ we may assume the following extra condition holds:

\begin{conv}\label{convention} Let $s\in C(L)$ with $L/K$ a finite extension and $\zeta=x(s)$. If $|\zeta|_w<R_w(\mathcal{Y})$ for some finite place $w\in \Sigma_{L,f}$, then there exists $1\leq t\leq l$ such that the sections $\tilde{s}$ and $\tilde{\xi}_t$ have the same image in $\tilde{C}(\F_{q(w)})$.
\end{conv}

With this convention in mind we make the following definition.

\begin{defn}\label{vadicproxdefn}
	Let $s\in C(L)$ and $w\in\Sigma_{L}$. We say that $s$ is \textbf{$w$-adically close to $s_0$} if $|x(s)|_w<\min\{ 1, R_w(\mathcal{Y})\}$.
\end{defn}

	\subsection{Strongly exceptional points}\label{section:strongexception}

In \cite{andre1989g} Andr\'e defines the notion of ``exceptional fibers'' of a one-parameter family $f:X\rightarrow S$ of abelian varieties as those fibers over points $s\in S(\C)$ for whom the corresponding endomorphism algebra $D_s$ is strictly larger than that of the generic fiber. He then looks at such fibers for which the further condition
\begin{equation}
	D_s\not\hookrightarrow M_g(\Q)
\end{equation}holds. Under the assumption that the abelian scheme has completely multiplicative reduction at $s_0$ he establishes the height bounds we are trying to establish in our context.

In analogy with Andr\'e's definition we will need, as eluded to in the introduction, certain conditions at the fibers reflecting the fact that we are no longer in the case where the reduction at $s_0$ is multiplicative. Given that our conditions are more complex, than Andr\'e's much more elegant analogue, we have gathered them in this section in the form of a definition.\\

Before we give the definition we introduce some notation concerning the data of the connected N\'eron model of the semiabelian variety $X'_0$. We consider the Chevalley decomposition of $X'_0$ as in \eqref{eq:chevalleysettingsection4} and write $B$ for the abelian part of $X'_0$, whose connected N\'eron model we denote by $\mathcal{B}$.

For any finite place $v\in\Sigma_{K,f}$ we write $h_v$ for the rank of the toric part of the fiber at $v$ of $\mathcal{B}$. We also write $B(v)$ for the abelian part of the same fiber. Finally, we assume that $B(v)$ is isogenous to $B_1(v)\times\cdots\times B_{r(v)}(v)$ where the $B_i(v)$ are the isotypic components of $B(v)$.

\begin{defn}\label{defnstrexc}Let $f:X\rightarrow S$, $s_0$, $g$, $h$, $B$, and $T$ be as in \Cref{section:settingnontrivial}. Let $s\in S(L)$, with $L/K$ finite, and assume that $X_s$ is isogenous to $X_1^{n_1}\times\ldots\times X_m^{n_m}$ where $X_i$ are simple non-isogenous abelian varieties and let $D_j:=\End^{0}(X_j)$.
	
The point $s$ will be called \textbf{strongly exceptional for the family} $f':X'\rightarrow S'$ if the following hold
	
\begin{enumerate}[label=(\roman*)]
	\item \underline{either} of the following two conditions are satisfied 
	\begin{enumerate}[label=(\alph*),ref=\theenumi(\alph*)]
		\item\label{strong1}there exists $j$ such that 
		\begin{center}
			$h> \frac{2\dim X_j}{[Z(D_j):\Q]}$ for some $j$, or
		\end{center}

		\item\label{strong2} $h\geq \min\{ 2  \frac{\dim X_j}{[Z(D_j):\Q] } : j \text{ such that }  D_j  \text{ is of type IV } \}.$

	\end{enumerate}

	\item\label{strong3}For all finite places $v\in \Sigma_{L,f}$ there exists $1\leq j\leq m$ such that there exists no embedding $M_{n_j}(D_j)\hookrightarrow M_{h+h_v}(\Q)$ and there exists no embedding $M_{n_j}(D_j)\hookrightarrow \End^{0}(B_i(v))$ for any of the $i$.
\end{enumerate}

\end{defn}

\begin{rmk}1. Given $A$ an abelian variety over some field $k$, we denote by $\End^{0}(A)$ the algebra $\End(A\times_k \bar{k})\otimes_{\Z}\Q$.\\
	
	2. As mentioned in the introduction \ref{strong1} and \ref{strong2} are based on the conditions outlined in $\S 8.4$  of \cite{papas2022height}, see in particular $\S8.4.1$ of loc. cit.. These are needed to define relations among the (archimedean) values of G-functions at $\xi=x(s)$, under the assumption that the point $s$ is archimedeanly close to the point $s_0$ for at least one archimedean place.\\
	
	3. On the other hand \ref{strong3} is needed to guarantee that the strongly exceptional points are not $v$-adically close to the point $s_0$ for any finite place $v$. We return to this in \Cref{section:globality}.\\
	
	4. We provide examples of families of strongly exceptional points in \Cref{section:examples}.\\
	
\end{rmk}

	\subsection{Relations at archimedean places}

We consider a fixed $f':X'\rightarrow S'$ as in \ref{section:settingnontrivial} and $f'_C:X_C'\rightarrow C'$ the semiabelian scheme associated to a good cover of $S'$. We also fix its associated family of G-functions $\mathcal{Y}$ and the subfamilies $\mathcal{Y}_t$, $1\leq t\leq l$, associated to each of the roots of the uniformizing parameter $x$ of the good cover $C_4$ of $S'$.

We introduce a final piece of notation, we will denote the elements of the subfamily $\mathcal{Y}_t$ by $y^{(t)}_{i,j}(x)$, where as usual $1\leq i\leq 2g$ and $1\leq j\leq h$ keep track of the entry in the relative period matrix that this power series corresponds to.\\

The construction in Lemma $8.3$ and the discussion in $\S 8.4.1$ of \cite{papas2022height} give us the following: 

\begin{lemma}\label{relarchplace1}Let $s\in C(L)$ be a strongly exceptional point of the G-admissible abelian scheme $f_C:X_C\rightarrow C$. Assume that $s$ is $w$-adically close to $s_0$ for some $w\in \Sigma_{L,\infty}$ with $w|v$. Then, there exists a root $\xi_t$ of the parameter $x$, where $1\leq t\leq l$, and a homogeneous polynomial \begin{center}
		$q_w\in \bar{\Q}[x_{i,j};1\leq i\leq 2g, 1\leq j\leq h]$,
	\end{center}with $\deg q_w\leq 2$, such that $q_w(i_w(y^{(t)}_{i,j})(x(s)))=0$.
\end{lemma}

\begin{rmk}In fact the polynomial $q_w$ will have coefficients in the extension $\hat{F}_s$ that depends on $L$ and the center of the algebra $\End^0(X_s)$. For more on this see \cite{andre1989g} or \cite{papas2022height}.
	
	Also from construction of these polynomials we know that they are not the specialization at $s$ of the polynomials defining the Riemann relations among the relative periods $y^{(t)}_{i,j}(x)$. For a description of these see \cite{papas2022height} and especially Lemma $7.5$.
\end{rmk}

\begin{proof}Since $s$ is $w$-adically close to $s_0$, by the properties of the good cover $C_4$ of $C'$ we get that there exists some root $\xi_t$ of $x$  and an analytic disk $\Delta_{w,R}\subset (C'\times_{w}\C)^{an}$ centered at $\xi_t$ with $R\leq R_w(\mathcal{Y})$ such that $s\in \Delta_{w,R}$. We then have that the values $i_w(y^{(t)}_{i,j})(x(s))$ appear as the entries of the first $h$ columns of a period matrix at $s$. This follows from \Cref{abeliangfunctions} and \Cref{abelianchangeofplace}.
	
	The relations among the periods in question follow from the definition of strongly exceptional points and the aforementioned constructions in loc. cit.
\end{proof}

\begin{defn}\label{defofrelations}
	Let $s\in C(L)$ be a strongly exceptional point of $f_C:X_C\rightarrow C$ and let $\zeta:=x(s)$. Then we define 
	\begin{equation}\label{eq:actualabnontrivial}
		q_s:=\prod_{\underset{|\zeta|_v<\min\{1, R_v(\mathcal{Y}) \}}{v\in \Sigma_{L,\infty }}}^{} q_v,
	\end{equation}where $q_v$ is the polynomial obtained in \Cref{relarchplace1}.
\end{defn}

\begin{rmks}1. For each archimedean place $v\in \Sigma_{L}$ for which $|\zeta|_v<\min\{1,R_v(\mathcal{Y})\}$ we have by construction of the $q_v$ in \Cref{relarchplace1}, that for the corresponding embedding $\iota_v:L\hookrightarrow \C$ \begin{equation}
		\iota_v(q_v)(\iota_v(y^{(t)}_{i,j})(\zeta))=0,
	\end{equation}where we note that the index $t$, indicating which root of $x$ the point $s$ is $v$-adically close to, depends on the place $v$! In particular we get that in the terminology of Ch. $VII$, $\S5$ of \cite{andre1989g}, the polynomial $q_s$ defines a relation among the values of the G-functions on the family $\mathcal{Y}$ at the point $\zeta=x(s)$ that hold $v$-adically for all archimedean places $v\in \Sigma_{L,\infty}$ with $|\zeta|_v<\min\{1, R_v(\mathcal{Y}) \}$.\\

2. What is missing in order for us to say that the relation induced from $q_s$ among the values of the aforementioned G-functions at $\xi$ are in fact ``global'', in the notation of loc. cit., is that they hold $v$-adically for all non-archimedean places $v\in \Sigma_{L,f}$ for which $s$ is $v$-adically close to the point of degeneration $s_0$. In the next section, following Andr\'e's argument we establish that strongly exceptional points cannot be $v$-adically close to $s_0$, thus establishing the globality of the above relations ``for free'' in this case.\\

3. We note that we have control over the degree of the polynomial $q_s$. Indeed, from construction of the polynomials $q_v$ we get that the degree of the polynomial $q_s$ is bounded above by $2 [L:\Q]$. See the end of the proof in $\S 14$ of \cite{papas2022height} for more details on this.
\end{rmks}

We finish this section by noting down the following lemma.

\begin{lemma}\label{nontrivcheck} Assume that for the G-admissible abelian scheme $f:X\rightarrow S$, assumed to be as in \Cref{section:settingnontrivial}, we also have that the generic special Mumford-Tate group associated to the variation of Hodge structures whose underlying local system is $R^1f^{an}_{*} \Q$ is the symplectic group $Sp(2g,\Q)$. Let also $C=c^{-1}(S)$ be the preimage of $S$ in some good cover of $S'$.
	
	Then for every strongly exceptional point $s\in C(\bar{K})$, that is archimedeanly close to the point $s_0$ for at least one archimedean place of its field of definition, the associated polynomial $q_s$ in \Cref{defofrelations} defines a non-trivial relation among the values of the associated family $\mathcal{Y}$ of G-functions at the point $\zeta=x(s)$.
\end{lemma}
\begin{proof}By the assumption that the point in question is $v$-adically close to $s_0$ for at least one archimedean place, the polynomial created does not define the ``$0=0$'' relation. 
	
Let us thus assume that the relation defined by $q_s$ is trivial, meaning that $q_s(\mathcal{Y})=0$ on the level of power series. Then by construction of $q_s$ we would have that one of the factors $q_v$ above is such that \begin{center}
	$q_v(\mathcal{Y}_t)=0$,
\end{center}for some $1\leq t\leq l$ depending on $v$, again on the level of power series. 

But in this case the subfamily we are working with, namely the family of power series $\mathcal{Y}_t$ comes from the degeneration of the semiabelian scheme $f_{C}':X'_C\rightarrow C'$ at the point $\xi_t$, which we remind the reader is just a root of the chosen uniformizing parameter $x$ of the good cover $C_4$ of $S'$. In this case the non-triviality follows by comparing the ideal defining the Zariski closure\footnote{See Lemma $7.5$ of \cite{papas2022height} for a description in this case.} of $\mathcal{Y}_t$ in $\mathbb{A}^{2g\times h}_{\bar{\Q}[[x]]}$  with the polynomials $q_v$. The non-triviality follows essentially by construction of the $q_v$, see $(42)$ and $(48)$ in the proof of Lemma $8.3$ in \cite{papas2022height}.\end{proof}


	\section{Globality}\label{section:globality}

We follow the notation set out in \Cref{section:settingnontrivial}, i.e. we consider $S'$ a smooth curve over a number field $K$, $s_0\in S'(K)$ a closed point, $f:X\rightarrow S$ an abelian scheme over $S:=S'\backslash\{s_0\}$ with $g$-dimensional fibers. As before we assume that the fiber of the connected N\'eron model $X'\rightarrow S'$ of $X$ over $S'$ has semiabelian fiber $X'_0$ whose Chevalley decomposition \begin{equation}
	0\rightarrow T\rightarrow X'_0\rightarrow B\rightarrow 0,
\end{equation}is such that $T$ is an $h$-dimensional split $K$-torus with $h\geq 2$ and $B$ is a $(g-h)$-dimensional abelian variety defined over $K$ with semiabelian reduction.

As in \Cref{section:settingnontrivial} we also fix a triple $(C_4,x,c)$ that keeps track of a good cover of $S'$. We also let $\tilde{C}$ be a regular projective model of the curve $C'$ over the ring $\mathcal{O}_K$ for which \Cref{convention} holds.

We establish the following analogue of Lemma $3.4$ in \cite{andre1989g}. Our proof follows in parallel the proof in loc.cit.

\begin{lemma}\label{gabberslemma}Assume that the point $s\in C(L)$ is strongly exceptional, where $L/K$ is a finite extension. Then there is no finite place $v\in \Sigma_{L,f}$ such that $s$ is $v$-adically close to $s_0$.
\end{lemma}

\begin{proof}Assume the opposite holds and let $v\in \Sigma_{L,f}$ be a finite place for which $s$ is $v$-adically close to $s_0$. Consider the sections $\tilde{s}$ and $\tilde{\xi}$, induced from the points $s$ and $\xi$ respectively, of the arithmetic pencil $\tilde{C}\times_{\spec\mathcal{O}_K}\spec{\mathcal{O}_L}$, where $\xi$ is some fixed root of the uniformizing parameter $x$ which is $v$-adically close to $s$ in the rigid-analytic sense.
By \Cref{convention} we know then that these sections will have the same image in $\tilde{C}(\mathbb{F}_{q(v)})$, where $q(v)$ is the cardinality of the residue field of $\mathcal{O}_{L_v}$.

We follow the main argument in \cite{andre1989g}. First of all, by spreading out there exists an open dense $\tilde{U}\subset \tilde{C}$ with $C\subset \tilde{U}$ such that the abelian scheme $f_C:X_C\rightarrow C$ descends to an abelian scheme $f_{\tilde{U}}:\tilde{X}_{\tilde{U}}\rightarrow \tilde{U}$.

By Gabber's Lemma, see $1.4$ and $4.10$ in \cite{delignegabber}, we have that there exists a proper surjective morphism $\psi: \tilde{T}\rightarrow \tilde{C}$ such that the abelian scheme $g_{\tilde{U}}:\psi^{-1}\tilde{X}_{\tilde{U}}:= \tilde{X}_{\tilde{U}}\times_{\tilde{U}} \psi^{-1}\tilde{U} \rightarrow \psi^{-1}\tilde{U}$ extends to a semiabelian scheme $\tilde{f}:\tilde{X}\rightarrow \tilde{T}$. By $1.7.b$ of \cite{delignegabber} the morphism $\psi$ may be furthermore assumed to be projective and the scheme $\tilde{T}$ to be integral and normal.              

Theorem $1.2$ in \cite{delignegabber} then shows that the extension $\tilde{X}\rightarrow \tilde{T}$ is unique and by the same result, as mentioned in \cite{andre1989g}, we have that 
$\tilde{X}\times_{\tilde{T}} \psi^{-1}(s)$ is the connected N\'eron model of $X_s$ over $\mathcal{O}_K$. Similarly, $\tilde{X}\times_{\tilde{T}} \psi^{-1}(s)$ is the connected N\'eron model of $X_s$ and $\tilde{X}\times_{\tilde{T}} \psi^{-1}(\xi)$ is the connected (lft) N\'eron model of the semiabelian variety $X'_{\xi}\simeq X'_{0}$. 

By our earlier assumptions on $T$ and $B$, see also \cite{halleniccomps} Proposition $4.1$, we get that $\tilde{X}\times_{\tilde{T}} \psi^{-1}(s_0)$ is a semiabelian scheme over $\mathcal{O}_K$. Since $\tilde{s}$ and $\tilde{s}_0$ have the same image in $\tilde{S}(\mathbb{F}_{q(v)})$, we get that the fibers of $\tilde{X}\times_{\tilde{T}} \psi^{-1}(s_0)$ and $\tilde{X}\times_{\tilde{T}}\psi^{-1}(s)$ at the place $v$ coincide.

From standard properties of N\'eron models we then have a natural morphism of $\Q$-algebras, induced by the N\'eron mapping property and specialization at $v$,\begin{equation}\label{eq:thereduction}
	\pi_v:\End^{0}(X_s)=\End^{0}(\tilde{X}\times_{\tilde{T}}\psi^{-1}(s))\rightarrow \End^{0}(X'_0(v)),
\end{equation}where $X'_0(v):=\tilde{X}\times_{\tilde{T}} \psi^{-1}(s_0)\times_{\spec{\mathcal{O}_L}}\spec\mathbb{F}_{q(v)}$ is the fiber of the above connected N\'eron model at the place $v$.

Assume that the Chevalley decomposition of $X'_0(v)$ is given by \begin{equation}\label{eq:chevalleyatprime}
	0\rightarrow T(v)\xrightarrow{i_v} X'_0(v)\xrightarrow{j_v} A(v)\rightarrow 0,
\end{equation}and let $h_v^{0}$ be the dimension of the torus $T(v)$.

Let $h_v$ be the toric rank at $v$ of the fiber of the connected N\'eron model of the abelian variety $B$ over the place $v$. Then by functoriality properties of N\'eron models we have that \begin{center}
		$h^{0}_v=h+h_v.$
	\end{center}
By functoriality properties of N\'eron models we also have that $A(v)$ is the abelian part of the fiber of the connected N\'eron model of the abelian variety $B$ over the place $v$.

From basic properties of algebraic groups, we have that the there exists an embedding of $\Q$-algebras
\begin{equation}\label{eq:semiabelianendomo}
	   \End^{0} (X'_0(v)) \hookrightarrow M_{h^{0}_v}(\Q)\oplus \End^{0} (A(v)).
\end{equation}

From \eqref{eq:thereduction} and \eqref{eq:semiabelianendomo} we get a homomorphism of $\Q$-algebras \begin{equation}\label{eq:punchline1}
	\End^{0}(X_s)\rightarrow M_{h^{0}_v}(\Q)\oplus \End^{0}(A(v)).
\end{equation}

On the other hand, we have that \begin{center}
	$\End^{0}(A(v))=\Bigsum{i=1}{r(v)} \End^{0}(B_i(v))$,
\end{center}where $B_i(v)$ are as in \Cref{section:strongexception}.

\begin{claim}\label{claimreductioninjective}The $\Q$-algebra homomorphism $\pi_v$ in \eqref{eq:thereduction} is injective.
\end{claim}

Assuming \Cref{claimreductioninjective} for the moment we complete the proof of \Cref{gabberslemma}.\\

Write $X_s\sim X_1^{n_1}\times \ldots \times X_m^{n_m}$ where the $X_i$ are mutually non-isogenous simple abelian varieties. Then $\End^{0}(X_s)= \Bigsum{i=1}{m} M_{n_i}(D_i)$ with $D_i:=\End^{0}(X_i)$. For each $1\leq i\leq m$ consider the composition of the injection $M_{n_i}(D_i)\hookrightarrow \End^{0}(X_s) $ with \eqref{eq:thereduction} and \eqref{eq:punchline1}. We then get an injective, thanks to \Cref{claimreductioninjective}, homomorphism of $\Q$-algebras\begin{equation}\label{eq:punchline2}
M_{n_i}(D_i)\hookrightarrow M_{h^{0}_v}(\Q)\oplus \End^{0}(A(v)).
\end{equation}

From simplicity of $M_{n_i}(D_i)$ we get that, for each $i$, either $M_{n_i}(D_i)\hookrightarrow M_{h^{0}_{v}}(\Q)$ or $M_{n_i}(D_i)\hookrightarrow\Bigsum{j=1}{r(v)} \End^{0}(B_j(v))$. The last embedding would once again imply by simplicity of $M_{n_i}(D_i)$ that $M_{n_i}(D_i)\hookrightarrow  \End^{0}(B_j(v))$ for some $j$.

Thus for each $i$ we have that either  $M_{n_i}(D_i)\hookrightarrow M_{h^{0}_{v}}(\Q)$ or $M_{n_i}(D_i)\hookrightarrow  \End^{0}(B_j(v))$ for some $j$.
This contradicts condition \ref{strong3} of strong exceptionality of the point $s$, see \Cref{defnstrexc}.
\end{proof}

\begin{proof}[Proof of \Cref{claimreductioninjective}] Let $l\neq p=p(v)$ be a prime, where $p(v)\in \N$ is the characteristic of the residue field of $L_v$. We also let $X_s[l^{\infty}]:=\cup_{n=1}^{\infty} X_s[l^n](\bar{L})$. Note that for $N\in \{ l^n:n\in \N\}$, we have that $X_s[N](\bar{L})\hookrightarrow X_s(v)$. Furthermore, for any $f\in \End^{0}(X_s)$ we trivially have that $f(P)\in X_s[l^{\infty}]$ for all $P\in X_s[l^{\infty}]$.
	
	Assume that $f\in \End^{0} (X_s)$ is such that $\pi_v(f)=0$.  Then we have that the reduction of  the set $f(X_s[l^{\infty}])$ at the place $v$ is trivial.  This implies that $f(X_s[l^{\infty}])=\{0\}$. But the set $X_s[l^{\infty}]$ is Zariski dense in $X_s$, see Theorem $5.3$ in \cite{edivdgmoo}, and hence $f=0$.
	
\end{proof}
	
	\section{Height bounds and some examples of strongly exceptional points}\label{section:everythingtogether}

In this section we put everything together. We record the analogue of Andr\'e's height bounds obtained in our case. We then discuss some examples that will hopefully elucidate the notion of strong exceptionality of a point $s$ of some abelian scheme $f:X\rightarrow S$. 

\subsection{Height bounds}

As usual we let $K$ be a number field and $S'$ a smooth connected curve over $K$, fix $s_0\in S(K)$ a closed point, and let $S=S'\backslash\{s_0\}$ be its complement. 

\begin{theorem}\label{heightboundsabelian}Let $f:X\rightarrow S$ be an abelian scheme with $g$-dimensional fibers, $g\geq 2$. We assume that the connected N\'eron model of $X$ over $S'$ is such that its fiber over $s_0$ has toric rank $h\geq 2$ and that the image of the morphism $S\rightarrow \mathcal{A}_g$ induced by $f$ is Hodge generic in $\mathcal{A}_g$.

Finally, we consider the set \begin{center}
	$\Sha(S):=\{s\in S(\overline{\Q}) : s \text{ is strongly exceptional for }f:X\rightarrow S     \}$.
\end{center}Then, there exist positive constants $c_1$, $c_2$ such that $h(P)\leq c_1 [K(P):K]^{c_2}$ for all $P\in \Sha(S)$.
\end{theorem}

This result will follow from the analogous result where $f':X'\rightarrow S'$ is replaced by the semiabelian scheme $f'_C:X_C'\rightarrow C'$ induced from a good cover of $S'$, see \Cref{goodcoversdef}.

\begin{theorem}\label{heightboundgoodcover}Let $f_C':X_C'\rightarrow C'$ be the semiabelian scheme associated to a good cover of the curve $S'$, where $S'$ is as in \Cref{heightboundsabelian}. 
	Consider the set \begin{center}
		$\Sha(C):=\{s\in C(\overline{\Q}) : s \text{ is strongly exceptional for }f_C:X_C\rightarrow C     \}$.
	\end{center}Then, there exist positive constants $c_3$, $c_4$ such that $h(P)\leq c_3[K(P):K]^{c_4}$ for all $P\in \Sha(C)$.
\end{theorem}
\begin{proof}Consider the set $\mathcal{Y}$ of G-functions associated to our cover, as in \Cref{associatedgfunctions}.
	
Let us fix from now on the point $s\in \Sha(C)$ and write $L:=K(s)$. From \Cref{gabberslemma} we know that for all $v\in \Sigma_{L,f}$ we will have that $|x(s)|_v\geq \min\{1,R_v(\mathcal{Y})\}$.  
		
Consider the set $\Sigma_\infty(s):=\{ v\in \Sigma_{L,\infty}: |x(s)|_v<\min\{1,R_v(\mathcal{Y}) \}\}$. If this set is empty then the argument laid out in the proof of Theorem $1.1$ of \cite{papas2022height} shows that the height $h(x(s))$ is bounded by an absolute constant independent of $s$.

If, on the other hand, the set $\Sigma_\infty(s)$ is non-empty, then \Cref{nontrivcheck} gives us non-trivial relations among the values of the G-functions of the family $\mathcal{Y}$ evaluated at $\zeta:=x(s)$. These relations are global, in the notation of Ch. $VII$, $\S5$ of \cite{andre1989g}, by virtue of \Cref{gabberslemma}. Hence Andr\'e's so called ``Hasse principle for values of G-functions'', see Theorem $5.2$ in Chapter $VII$ of \cite{andre1989g}, applies in our case.

The aforementioned result shows that $h(x(s))\leq c_5 \delta^{c_6}$ where $\delta$ is the degree of the polynomial defined in \Cref{defofrelations}. These degrees are bounded by construction of these polynomials by $2[\hat{F}_s:\Q]$ where $\hat{F}_s$ is a finite extension\footnote{For the definition of $\hat{F}_s$ see the discussion at the end of section $\S8.1$ of \cite{papas2022height}.} of $L$, in the notation of \cite{andre1989g} this is the field $\hat{K}$ on page $202$. Results of Silverberg, see \cite{silverberg}, then imply that $[\hat{F}_s:\Q]\leq c_{7}[L:K]$ finishing the proof.
\end{proof}

\begin{proof}[Proof of \Cref{heightboundsabelian}] Let $s\in \Sha(S)$ be a strongly exceptional point of $f:X\rightarrow S$. Then by construction of the good cover $C_4$, see \Cref{section:goodcovers}, we know that there exists $s_1\in C(L)$ such that $(X_C)_{\tilde{s_1} }\simeq X_s$. In particular, since the fibers $(X'_C)_{\xi}$, where here $\xi$ ranges over the roots of $x$, are isomorphic to $X'_{s_0}$, it is easy to see that $s_1\in\Sha(C)$. 
		
Therefore the result follows from \Cref{heightboundgoodcover}, since for any Weil height $h_S$ on $S'$ we will have that there exists a positive constant $c_5$ such that \begin{center}
	$h_S(s)\leq h_C(s_1)+c_5$,
\end{center}where $h_C$ is the Weil height induced from $h_S$ on $C'$ via the covering morphism\footnote{See \Cref{section:goodcovers} for the notation used here.} $c:C_4\rightarrow \bar{S}'$.
	
\end{proof}

	\subsection{Examples}\label{section:examples}

We present here some notable examples of strongly exceptional points that are covered by \Cref{heightboundsabelian}. We follow here the general notation set out in \Cref{section:settingnontrivial}. In particular, as usual we denote the fiber over $s_0$ of the connected N\'eron model of $X$ by $X'_0$ and let $T$ be its toric part and $B$ be its abelian part.

\subsubsection*{Abelian part is 1-dimensional with everywhere good reduction}

The conditions, especially \ref{strong3}, of \Cref{defnstrexc} are especially easy to check if one imposes strong assumptions about the reduction of the abelian part $B$ of $X_0$. In this subsection we assume that $B$ is an \textbf{elliptic curve that has everywhere good reduction\footnote{By \cite{serretate} all CM elliptic curves over a number field satisfy this relation after base change by a finite extension of our ground field $K$.}}. Note that in this case we will have $h=h^0_v=g-1$ for all finite places $v$.

Let us assume that $X_s\sim X_1^{n_1}$ with $X_1$ a simple abelian variety. Write $D_1:=\End^{0}(X_1)$, so that we have $\End^{0}(X_s)=M_{n_1}(D_1)$, and let $e:=[Z(D_1):\Q]$. We are particularly interested, with a view towards the results in \cite{daworr3}, in the case where $n_1=1$.

Note that in this case, i.e. $n_1=1$ and $X_s$ is simple, \ref{strong1} is equivalent to $e> 2+\frac{2}{g-1}$. Upon assuming $g\geq 4$, which is where we deviate from previous known height bounds, see \cite{andre1989g} and \cite{daworr}, this is equivalent to $e\geq 3$.\\

From classical results, see the table on page $187$ of \cite{mumfordabelian}, we must have that $e|g$ when $D_1$ is of types $I-III$. Therefore we have that $D_1\not\hookrightarrow M_h(\Q)$ since $e\not| h$. All that is left to establish \ref{strong3} is to check that $D_1\not\hookrightarrow \End^{0}(B(v))$ for all finite places $v$. Upon assuming that $e\geq 3$ this becomes trivial by the description of the endomorphism algebra of an elliptic curve over a finite field, see \cite{silverman1}. 

\begin{lemma}\label{example1}Let $f:X\rightarrow S$ be a G-admissible abelian scheme with fibers of dimension $g$. We assume that the abelian part $B$ of $X'_0$ is as above, i.e. an elliptic curve with everywhere good reduction, and that $g\geq 5$. Then any point $s\in S(\overline{\Q})$ for which $X_s=X_1$ is simple of type $I-III$ in Albert's classification and $e=[Z(D_1):\Q]\geq 3$ is strongly exceptional.
\end{lemma}

Another interesting phenomenon that appears in this case, upon assuming that $g\geq 6$ is even, is that points $s$ whose fibers $X_s$ are isogenous to $X_1^2$ with $X_1$ a simple abelian variety can be strongly exceptional. This case is not covered in Y. Andr\'e's original result, or the subsequent result of Daw and Orr in \cite{daworr} that covers the case where $g$ is even, under the assumption of having completely multiplicative reduction at $s_0$.

\begin{lemma}\label{example2}Let $f:X\rightarrow S$ be a G-admissible abelian scheme with fibers of dimension $g$. We assume that the abelian part $B$ of $X'_0$ is as above, i.e. an elliptic curve with everywhere good reduction, and that $g$ is even with $g\geq 6$. Then any point $s\in S(\overline{\Q})$ with $X_s\sim X_1^{2}$, where $X_1$ is simple of type $I-III$ in Albert's classification, with $e= [Z(D_1):\Q]\geq 3$ is strongly exceptional.
\end{lemma}

\begin{proof}From the fact that $e\geq 3$ we immediately get that \ref{strong1} and \ref{strong3} is satisfied. Here we are using the assumption that $B$ has everywhere good reduction and that $e|g/2$, and hence $M_2(D_1)\not\hookrightarrow M_{g-1}(\Q)$.
\end{proof}

\begin{rmks}1. In Andr\'e's classical result the nature of the reduction, i.e. completely multiplicative, does not allow points such as those considered in \Cref{example2}. This issue remains in the analogous result of Daw and Orr, see Theorem $8.1$ of \cite{daworr}. 
	
	Indeed, assume that the reduction is completely multiplicative\footnote{In other words, $h=g$ in our notation and $X'_0$ is a $g$-dimensional torus} and let $s\in S(\bar{\Q})$ is as in \Cref{example2} with $D_1$ of type $I$. In that case we would need to have that there is no embedding \begin{center}
		$\End^{0}(X_s)\hookrightarrow M_g(\Q)$.
	\end{center}But since $X_s$ is isogenous to $X_1^2$ whose algebra $D_1$ is of type $I$, we must have $e|g/2$ which implies  $D_1\hookrightarrow M_{g/2}(\Q)$. On the other hand $\End^0(X_s)= M_2(D_1)$ and hence an embedding into $M_g(\Q)$ exists making it impossible to continue.\\
	
	2. We note that Daw and Orr in \cite{daworr4} deal with such isogenies in the case where the abelian scheme is a product of elliptic schemes over $S$, again assuming completely multiplicative reduction over a point.
	
	To achieve this the critical new ingredient is the introduction of relations among values of G-functions at finite places.\\

	3. We note that the above corollaries can easily be replicated with $B$ no longer an elliptic curve but a higher dimensional abelian variety that has everywhere good reduction. To extract aesthetically pleasing examples one only has to allow the dimension $g$ of the fibers to be ``big enough''.
\end{rmks}
	
	\section{Some cases of the Zilber-Pink Conjecture}\label{section:caseszp}

In \cite{daworr3}, generalizing their work in \cite{daworr, daworr2},  C. Daw and M. Orr consider curves $S\subset \mathcal{A}_g$ with $g\geq 3$. They consider intersections of such curves with the set of points $\Sigma \subset \mathcal{A}_g(\C)$ that correspond to simple abelian varieties with endomorphism algebra of type $I$ or $II$ in Albert's classification.

The Zilber-Pink conjecture predicts that the number of these intersections is finite assuming the curve $S$ in question is Hodge generic. C. Daw and M. Orr reduce the Zilber-Pink conjecture in this setting to a Large Galois orbits hypothesis, see Conjecture $8.2$ in \cite{daworr3}. Following work in \cite{daworr}, namely the proof of Proposition $9.3$ in loc.cit., the Large Galois orbit hypothesis can be reduced to height bounds of the type we get in \Cref{heightboundsabelian}.

We start by establishing Conjecture $1.5$ of \cite{daworr3} for strongly exceptional points in curves such as those appearing in \Cref{heightboundsabelian}. We close off our exposition by giving some more concrete examples of Zilber-Pink type statements that follow from \Cref{largego}.

\subsection{Large Galois Orbits}
	
	\begin{cor}[LGO]\label{largego} Let $f:X\rightarrow S$ be a G-admissible\footnote{See \Cref{section:notation} for our conventions on G-admissible abelian schemes.} abelian scheme defined over the smooth irreducible curve $S'$, with $S=S'\backslash \{s_0\}$, with everything defined over a number field $K$.
	
	 Let $\Sha(S)$ be as in \Cref{heightboundsabelian} and consider the sets 
	\begin{center}
		$\Sha_0(S):=\{s\in S(\C) :X_s \text{ is simple } \}$, and\end{center}
	
	\begin{center}		
		$\Sha_1(S)=\Sha_0(S)\cap \Sha (S)$.
	\end{center}

Then there exist positive constants $c_5$, $c_6$ such that for all $s\in 	\Sha_1(S)$ we have \begin{center}
	$\# \aut(\C/K)\cdot s\geq c_5 |\disc(\End(X_s))|^{c_6}$
\end{center}
\end{cor}

\begin{proof} The general strategy to reduce LGO to the height bounds derived from the G-functions appears in the work of Daw and Orr, \cite{daworr, daworr2,daworr3}. We chose to present a summary of the argument for reasons of completeness of our exposition.
	
	We note that we first need to base change everything by a finite extension $K'/K$, so that \Cref{conventionsemistable} holds.
	
From \Cref{heightboundsabelian} we get that there exist positive constants $c_1, c_2$ such that for all $s\in \Sha_1(S)$, after choosing a Weil height $h$ on $S'$, we have \begin{equation}\label{eq:weilvsdeg}
	h(s)\leq c_1 [K(s):K]^{c_2}.
\end{equation}

The Weil height and the stable Faltings height $h_F$ have the following\footnote{This is proven in Lemma $3$ of \cite{faltings}.} comparison\begin{equation}\label{eq:weilvsfaltings}
	|h_F(X_s)-h(s)|=O(\log (h(s))).
\end{equation}Combining those one gets the existence of positive constants $c_8$, $c_9$ such that \begin{equation}\label{eq:faltvsdegree}
h_F(X_s)\leq c_7[K(s):K]^{c_8},
\end{equation}for all $s\in \Sha_1(S)$.

The final piece that we need is Theorem $1$ of \cite{maswus}. This result, paired with Lemma $5.6$ of \cite{daworr2}, gives that \begin{equation}\label{eq:maswusest}
	|\disc(\End (X_s))|\leq c_{9} \max([K(s):K], \delta, h_F(X_s))^{c_{10}},
\end{equation}where $c_{9}$, $c_{10}$ are positive constants depending only on $g$, the dimension of the abelian scheme $X_s$, and $\delta$ is the degree of a polarization on $X_s$. For more on this see the exposition in \cite{daworr2} page $55$. 

Finally, from Theorem $1.1$ of \cite{poletiso} one can bound $\delta$ from above by a quantity, denoted $\kappa(X_s)$ in the notation of loc. cit., that only depends polynomially on $h_F(X_s)$ and $[K(s):K]$. Combining \eqref{eq:maswusest} with the aforementioned bound for $\delta$ and \eqref{eq:faltvsdegree}, finishes the proof of the bounds in question.
\end{proof}

\begin{rmk}
	One could consider in \Cref{largego} the intersection of $\Sha_0(S)$ with the subset $\Sha(S)'\subset S(\C)$ of strongly exceptional points, not ``necessarily'' in $\bar{\Q}$. In this case one gains nothing more since the points in  question will already belong in $\Sha_1(S)$, i.e. they will be points in $S(\bar{\Q})$. This follows because the points in $\Sha(S)'$ can be viewed as points of intersection of the image of the induced morphism $i:S\rightarrow \mathcal{A}_g$ with special subvarieties of $\mathcal{A}_g$. Viewed as subvarieties of $\mathcal{A}_g$ these are all defined over $\bar{\Q}$.
\end{rmk}

\subsection{Application to the Zilber-Pink Conjecture}\label{section:applicationtozp}

Before establishing \Cref{zpintro} we need the following analogue of Proposition $9.4$ of \cite{daworr}:

\begin{lemma}\label{lemma:auxioliarypullback}Let $Z\subset \mathcal{A}_g$ be an irreducible algebraic curve such that its Zariski closure in the Baily-Borel compactification intersects the $\mathcal{A}_{g-h}$-stratum of the boundary. Let $s^{*}$ be a fixed point in this intersection.
	
	Then there exists a smooth projective curve $\tilde{C}$, an open subset $C\subset \tilde{C}$, a point $s_0\in \tilde{C}\backslash C$, a finite surjective morphism $q:C\rightarrow Z$ and a semiabelian scheme $\tilde{f}:\tilde{X}\rightarrow \tilde{C}$ such that:\begin{enumerate}
		\item $X:=\tilde{X}|_{C}$ is an abelian scheme over $C$, 
		
		\item the map $C\rightarrow \mathcal{A}_g$ induced from $X\rightarrow C$ is the composition $C\xrightarrow{q} Z\hookrightarrow \mathcal{A}_g$,
		
		\item the fiber $\tilde{X}_{s_0}$ is a semiabelian variety with toric rank $h$,
		
		\item the abelian part of the semiabelian variety $\tilde{X}_{s_0}$ is isomorphic to the $(g-h)$-dimensional abelian variety whose classifying point is the above point $s^{*}\in \mathcal{A}_{g-h}$.
	\end{enumerate} 

We write $C':=\tilde{C}\backslash((\tilde{C}\backslash C)\backslash\{s_0\})$, i.e. $C':=C\cup \{s_0\}$, and $f':X'\rightarrow C'$ for the base change of $\tilde{f}$ along the open immersion $C'\rightarrow \tilde{C}$. We call $f':X'\rightarrow C'$ the \textbf{semiabelian scheme associated to the pair} $(Z,s^{*})$.
\end{lemma}

\begin{proof}We follow the proof on Proposition $9.4$ in \cite{daworr} and its general notation. So we denote by $\mathcal{A}_g^{*}$ the Baily-Borel compactification of $\mathcal{A}_g$ and by $Z^{*}$ the Zariski closure of the curve $Z$ in $\mathcal{A}_g^{*}$. We also let $s^{*}$ be the point of intersection of $Z^{*}$ with the $\mathcal{A}_{g-h}$-stratum of the boundary, whereas $s^{*}$ in loc. cit. is just the $0$-dimensional stratum of the fiber. We also let $\mathfrak{G}$, and $\mathfrak{G}_3$ be the semiabelian schemes over $\bar{\mathcal{A}}_g$ and  $\bar{\mathcal{A}}_{g,3}$ respectively as in loc. cit..
	
	The first part of the proof in loc. cit. passes through verbatim and gives an irreducible curve $V_3\subset \mathcal{A}_{g,3}$ whose Zariski closure $\bar{V}_3$ in some toroidal compactification $\bar{\mathcal{A}}_{g,3}$ of $\mathcal{A}_{g,3}$ has a point $\bar{s}_3$ which maps to $s^{*}$ under the map $\bar{\pi}\circ p$ of loc. cit..
	
	Likewise the second part of the proof in loc. cit. also passes through almost verbatim. In particular, we get a smooth projective curve $\tilde{V}'$ together with an open subset $\tilde{V}\subset\tilde{V}'$ and a semiabelian scheme $\mathfrak{U}'$ over the curve $\tilde{V}'$ as in loc. cit. and then consider $s_0\in \tilde{V}'$ that is a preimage of the point $\bar{s}_3$. 
	
	The compatibility of the map $\bar{\pi}$ in loc. cit. with the stratifications of $\mathcal{A}_{g}^{*}$ and $\bar{\mathcal{A}}_g$ implies that the fiber $(\mathfrak{G}_3)_{\bar{s}_3}$ is a semiabelian variety with toric rank $h$. 
	
	Thus taking $C:=\tilde{V}$, $\tilde{C}:=\tilde{V}'$ and $\tilde{X}:=\mathfrak{U}'$ parts $(1)-(3)$ of the Lemma follow.
	
	For part $(4)$ we note that, by Theorem $2.3.(4)$ of Ch. $V$ of \cite{falchai}, we have that $\mathfrak{G}_{q(\bar{s}_3)}$ is semiabelian with abelian part as in $(4)$. Thus the result follows by construction of the semiabelian scheme $\tilde{X}$.
	
\end{proof}

\subsubsection{Proof of \Cref{zpintro}}
We start with an obvious definition that follows from our work so far. We then establish \Cref{zpintro}, which is the most general Zilber-Pink type statement that we can establish.\\

Motivated by \Cref{lemma:auxioliarypullback}, we make the following:
\begin{defn}\label{stronglyexceptionalag}Let $Z\subset \mathcal{A}_g$ be an irreducible algebraic curve. Let $q:C\rightarrow Z$ and $f':X'\rightarrow C'$ be the curve and semiabelian scheme associated to the pair $(Z,s^{*})$ from \Cref{lemma:auxioliarypullback}, where $s^{*}$ is as in \Cref{lemma:auxioliarypullback}.
	
	We say that a point $s\in Z(\bar{\Q})$ is \textbf{strongly exceptional for $Z$} if there exists $\tilde {s}\in C(\bar{\Q})$ such that $q(\tilde{s})=s$ and $\tilde{s}$ is strongly exceptional for the semiabelian scheme $f':X'\rightarrow C'$.
\end{defn}

The general strategy of the proof of \Cref{zpintro} follows the discussion in \cite{daworr2} $\S 6.7$, we include a summary of the argument for the sake of completeness of our exposition.

\begin{proof}[Proof of \Cref{zpintro}]Consider $Z\subset \mathcal{A}_g$ as in the statement. We also let $\Sigma$ be the subset of $\mathcal{A}_g(\C)$ whose points correspond to abelian varieties that are simple and whose algebras of endomorphisms are of type $I$ or $II$ in Albert's classification and satisfy the conditions in \Cref{defnstrexc}.
	
	 Assume that $K_1$ is some number field of definition of $Z$. There exists a finite extension $K$ of $K_1$ such that the auxiliary curve $C'$, the point $s_0$, and the morphism $f':X'\rightarrow C'$ associated to $Z$ via \Cref{lemma:auxioliarypullback} are all defined over $K$. Note that the family $f:X\rightarrow C$ will be a G-admissible abelian scheme by virtue of the factorization of the covering $q$ in \Cref{lemma:auxioliarypullback}.
	
Consider $\Sha_0'(C)\subset \Sha_0(C)$, with notation as in \Cref{largego}, to be the subset of strongly exceptional points of the curve $C$ whose fibers $X_s$ are simple abelian varieties with endomorphism algebras of type $I$ or $II$ in Albert's classification. Then by \Cref{largego} we have that there exist positive constants $c_5$, $c_6$ such that \begin{equation}\label{eq:largegalois}
	\# \aut(\C/K)\cdot P\geq c_5 |\disc(\End(X_s))|^{c_6},
\end{equation}for all $P\in \Sha_0'(C)$.

Now we transfer this back to the original curve $Z$. Indeed, the finiteness of the morphism $q$ implies that if $s\in Z(\C)\cap \Sigma$, then we must have that $s\in Z(\bar{\Q})$, by virtue of being the point of intersection of a special subvariety of $\mathcal{A}_g$, which is thus defined over $\bar{\Q}$, with $Z$, itself defined over $\bar{\Q}$ by assumption.

Finiteness of the morphism $q$ of \Cref{lemma:auxioliarypullback} implies that there exists a point $\tilde{s}\in C(\bar{\Q})$ such that $\tilde{s}$ is strongly exceptional for the semiabelian scheme $f':X'\rightarrow C'$, by virtue of the factorization of the map $C\rightarrow \mathcal{A}_g$ via $Z\hookrightarrow \mathcal{A}_g$, and $q(\bar{s})=s$. This argument also shows that the definition of strongly exceptional points of $Z$ is independent of the cover $C$ chosen by \Cref{lemma:auxioliarypullback}.

By finiteness of the morphism $q$ we also get the existence of a positive constant $c_{11}$, independent of the points chosen, such that \begin{equation}\label{eq:morphgalois}
	\#\aut(\C/K_1)\cdot s\geq c_{11} \#\aut(\C/K)\cdot \tilde{s}.\end{equation}Combining \eqref{eq:largegalois} together with \eqref{eq:morphgalois} we get that there exist constants $c_{12}$, $c_{13}$ independent of $s$ such that\begin{equation}
	 \#\aut(\C/K_1)\cdot s\geq c_{12} |\disc\End (A_s)|^{c_{13}},
\end{equation}where $A_s$ is the abelian variety parameterized by the point $s$.

In other words, the curve $Z\subset \mathcal{A}_g$ satisfies Conjecture $1.5$ of \cite{daworr3}, with $\Sigma$ in loc. cit. replaced by the set $\Sigma$ we are considering. Then, the result follows from Theorem $1.3$ of \cite{daworr3}, which reduces the Zilber-Pink conjecture for intersections of the curve $Z$ with special subvarieties of simple PEL type\footnote{For a definition see section $1.C$ of \cite{daworr3}.} $I$ or $II$ to the aforementioned Conjecture of loc. cit., the so called ``large Galois orbits hypothesis'' in this setting.
\end{proof}

		\subsubsection{Examples of Zilber-Pink type statements}

We dedicate this last part to presenting some examples that follow from the much more general \Cref{zpintro}. In the process we see that the definition of ``strongly exceptional points'' is in fact far less restricting in some cases.\\

Inspired by Andr\'e's definition\footnote{See \Cref{section:strongexception} for more on this.} of exceptional fibers in $1$-parameter families of abelian varieties we introduce the following:
\begin{defn}Let $S\subset \mathcal{A}_g$ be a curve defined over some subfield $k$ of $\C$, and let $f:X\rightarrow S$ be the pullback of the universal abelian variety over $S$. A point $s\in S(\C)$ for which the fiber $X_s$ is exceptional and also a simple abelian variety will be called a \textbf{simple exceptional point of $S$}.
\end{defn}

Before we state our examples we introduce a useful bit of notation:

\begin{notation}Let $s\in \mathcal{A}_g(\C)$ we write $A_s$ for the associated $g$-dimensional abelian variety. If the abelian variety $A_s$ is \textbf{simple} we denote by $e_s$ the degree $[Z(D_s):\Q]$, where $D_s:=\End^{0}(A_s)$ is the algebra of endomorphisms of the abelian variety $A_s$.\end{notation}

\begin{cor}\label{zpprimes} Let $S\subset \mathcal{A}_g$, where $g\geq 5$ is an odd prime, be a Hodge generic irreducible algebraic curve defined over $\bar{\Q}$.
	
	Assume that the intersection of the Zariski closure of $S$ in the Baily-Borel compactification of $\mathcal{A}_g$ with some $\mathcal{A}_{g-h}$-stratum, where $3\leq h\leq g-1$, is a CM-point. Then there are only finitely many simple exceptional points of $S$ in $S(\bar{\Q})$ with $e_s\geq 3$.
\end{cor}
\begin{proof}We employ \Cref{lemma:auxioliarypullback} to shift focus from $S$ to a the semiabelian scheme $f':X'\rightarrow C'$ as in the aforementioned lemma. Note that, by part $(4)$ of \Cref{lemma:auxioliarypullback}, we have that the abelian part of $X'_0$, the fiber of $X'$ over the a preimage of the point of intersection of the Zariski closure $S^{*}$ of $S$ with the $\mathcal{A}_{g-h}$-stratum of the boundary of the Baily-Borel compactification $\mathcal{A}_g^{*}$, will be a CM abelian variety.
	
	Let $s\in C(\bar{\Q})$ be a preimage of a simple exceptional point of $S$. We use a case-by-case argument depending on the type of the simple fiber $X_s$ in Albert's classification.

	Since the dimension of the fibers $g$ is prime, by the table on page $187$ of \cite{mumfordabelian} we have that there are no fibers that are simple of type $II$ or $III$ in Albert's classification.
	
	If $s$ is a point whose corresponding fiber is simple of type $IV$ then, by the same table of loc. cit. and since the dimension $g$ is prime, we have that either $e=2g$, and hence $s$ is a CM-point, or that $e=2$. The case $e=2$ is ruled out by our assumptions, while the finiteness of CM-points follows from the Andr\'e-Oort conjecture, i.e. from the results of \cite{tsimermanag} in our case.
	
	We are thus reduced to points $s\in C(\bar{\Q})$ whose fiber is simple of type $I$. In this case we claim that the result follows from \Cref{zpintro}, in other words, we claim that in this case all simple type $I$ points are strongly exceptional.
	
	 Indeed, by the same table in \cite{mumfordabelian}, we have that the only possibility for the dimension of the algebra $D_s:=\End^0(X_s)$ is that $e_s=g$. Hence, since $h\geq 3$ by assumption, the point $s$ satisfies condition \ref{strong1} of \Cref{defnstrexc}. 

Since, by assumption, the abelian part of the semiabelian variety $X'_0$ is CM, and hence has potentially good reduction everywhere by \cite{serretate}, the quantities $h_v$ in condition \ref{strong3} of \Cref{defnstrexc} are automatically $h_v=0$ for all finite places $v\in \Sigma_{K,f}$. Hence we get, by the fact that $g$ is prime and $h\leq g-1$, that there is no embedding $D_s\hookrightarrow M_{h+h_v}(\Q)=M_{h}(\Q)$. 

Finally, we note that, since $g$ is prime, there cannot be any embedding $D_s\hookrightarrow \End^{0}(B_i(v))$, where $B_i(v)$ are as in \Cref{section:strongexception}. Hence, condition \ref{strong3} of \Cref{defnstrexc} is also satisfied and the result follows from \Cref{zpintro}.
\end{proof}

\begin{rmks} We note that we may replace the condition in \Cref{zpprimes} that the Zariski closure of the curve $S$ in the Baily-Borel compactification $\mathcal{A}_g^{*}$ intersects some stratum as above in a CM-point with the condition that it intersects such a fiber in a point with potentially good reduction everywhere.
\end{rmks}

Following the discussion in \Cref{section:examples} we may formulate, analogous to the example highlighted in \Cref{example1}, cases of the Zilber-Pink conjecture that follow directly from \Cref{zpintro}.

\begin{cor}\label{zpexample1}Let $S\subset \mathcal{A}_g$, where $g\geq 5$, be a Hodge generic irreducible algebraic curve defined over $\bar{\Q}$.
	
Assume that the intersection of the Zariski closure of $S$ in the Baily-Borel compactification of $\mathcal{A}_g$ with the $\mathcal{A}_{1}$-stratum is a CM-point. Then the set \begin{center}
	$\{ s\in S(\bar{\Q}):s $ is simple of type $I$ or $II$ and $e_s\geq 3 \}$
\end{center}is finite. 
\end{cor}
\begin{proof}The points of the set in question are strongly exceptional by \Cref{example1}. The result hence follows from \Cref{zpintro}.	
\end{proof}

\begin{rmks}1. As remarked earlier, see the end of \Cref{section:examples}, we can create examples, similar to the one above, with curves $S$ intersecting the boundary at some $A_{g-h}$ stratum with $g-h>1$. The difference is that we would have to choose larger $g$ and $h$ to get results of similar simplicity.\\
	
	2.  Let us assume that in the context of \Cref{zpexample1}, $g$ is also odd. Then, using similar arguments as in the proof of \Cref{zpprimes}, we can prove that there are only finitely many simple points of type $I-III$. Obviously in this case types $II$ and $III$ would be a non-factor in our analysis.\\
	
	3. Establishing the finiteness of simple type $III$ and $IV$ points that are strongly exceptional will in general require that the analogues of the results of \cite{daworr3}, used to establish Theorem $1.3$ of \cite{daworr3} for type $I$ and $II$ simple points, hold for simple type $IV$ points as well.\\
\end{rmks}
	



	\bibliographystyle{alpha}
	
	\bibliography{biblio}

\end{document}